\documentclass[english,journal]{IEEEtran}
\usepackage[T1]{fontenc}
\usepackage[utf8]{inputenc}
\usepackage{color}
\usepackage{verbatim}
\usepackage{mathrsfs}
\usepackage{mathtools}
\usepackage{url}
\usepackage{enumitem}
\usepackage{amsmath}
\usepackage{amsthm}
\usepackage{amssymb}
\usepackage{xargs}[2008/03/08]

\makeatletter

\providecommand{\tabularnewline}{\\}

  \theoremstyle{definition}
  \newtheorem{defn}{\protect\definitionname}
\theoremstyle{plain}
\newtheorem{thm}{\protect\theoremname}
  \theoremstyle{plain}
  \newtheorem{cor}{\protect\corollaryname}
  \theoremstyle{plain}
  \newtheorem{prop}{\protect\propositionname}
  \theoremstyle{remark}
  \newtheorem{rem}{\protect\remarkname}
  \theoremstyle{plain}
  \newtheorem{lyxalgorithm}{\protect\algorithmname}

\usepackage{color}
\usepackage{graphics} 
\usepackage{epsfig} 
\usepackage{amsfonts}
\usepackage{gensymb}
\usepackage{cite}
\allowdisplaybreaks

\makeatother

\makeatletter
\def\ps@headings{\def\@oddhead{\IEEEdoarxivheader{-1\oddsidemargin}\relax
\hbox{}\@IEEEheaderstyle\rightmark\hfil\thepage}\relax
\def\@evenhead{\IEEEdoarxivheader{-1\evensidemargin}\relax
\hbox{}\@IEEEheaderstyle\rightmark\hfil\thepage}\relax
\def\@oddfoot{\IEEEdoarxivfooter{-1\oddsidemargin}\hfil\hbox{}}\relax
\def\@evenfoot{\IEEEdoarxivfooter{-1\evensidemargin}\hfil\hbox{}}\relax}
\def\ps@IEEEtitlepagestyle{\ps@headings}

\def\IEEEarxivheadfootoffset{3pt}
\newdimen\IEEEheadtotopofpage
\newdimen\IEEEfoottobottomofpage
\newbox\@IEEEboxX

\def\IEEEarxivheader{}
\def\IEEEarxivfooter{}

\advance\IEEEheadtotopofpage by 1\topmargin
\advance\IEEEheadtotopofpage by 1in\relax

\IEEEfoottobottomofpage\paperheight
\advance\IEEEfoottobottomofpage by -1in\relax
\advance\IEEEfoottobottomofpage by -1\topmargin
\advance\IEEEfoottobottomofpage by -1\headheight
\advance\IEEEfoottobottomofpage by -1\headsep
\advance\IEEEfoottobottomofpage by -1\textheight
\advance\IEEEfoottobottomofpage by -1\footskip
\def\IEEEarxivheaderstyle{\normalfont\footnotesize}

\def\IEEEdoarxivheader#1{\@IEEEtrantmpdimenA\IEEEarxivheadfootoffset\relax
\@IEEEtrantmpdimenA -1\@IEEEtrantmpdimenA
\advance\@IEEEtrantmpdimenA by \IEEEheadtotopofpage
\settoheight{\@IEEEtrantmpdimenB}{\IEEEarxivheaderstyle HT}\relax
\advance\@IEEEtrantmpdimenA by -1\@IEEEtrantmpdimenB
\setbox\@IEEEboxX=\hbox{\relax
\raisebox{\@IEEEtrantmpdimenA}[0pt][0pt]{\parbox[t]{\textwidth}{\centering
\IEEEarxivheaderstyle\IEEEarxivheader}}}\relax
\wd\@IEEEboxX=0pt\relax
\ht\@IEEEboxX=0pt\relax
\dp\@IEEEboxX=0pt\relax
\box\@IEEEboxX\relax}

\def\IEEEarxivfooterstyle{\normalfont\footnotesize}

\def\IEEEdoarxivfooter#1{\@IEEEtrantmpdimenA\IEEEfoottobottomofpage\relax
\advance\@IEEEtrantmpdimenA by \IEEEarxivheadfootoffset\relax
\settodepth{\@IEEEtrantmpdimenB}{\IEEEarxivheaderstyle gjpqy}\relax
\advance\@IEEEtrantmpdimenA by 1\@IEEEtrantmpdimenB
\setbox\@IEEEboxX=\hbox{\hskip#1\hskip -1in\relax
\raisebox{\@IEEEtrantmpdimenA}[0pt][0pt]{\parbox[b]{\paperwidth}{\centering
\IEEEarxivfooterstyle\IEEEarxivfooter}}}\relax
\wd\@IEEEboxX=0pt\relax
\ht\@IEEEboxX=0pt\relax
\dp\@IEEEboxX=0pt\relax
\box\@IEEEboxX\relax}

\def\@IEEEheaderstyle{\normalfont\scriptsize}
\def\@IEEEfooterstyle{\normalfont\scriptsize}

\makeatother

\renewcommand{\IEEEarxivheadfootoffset}{3pt}

\renewcommand{\IEEEarxivheaderstyle}{\normalfont\footnotesize}
\renewcommand{\IEEEarxivfooterstyle}{\normalfont\footnotesize}

\renewcommand{\IEEEarxivheader}{This is the author's version of an article
that has been accepted to IEEE Transactions on Automatic Control. Changes were made to this version
by the publisher prior to publication.\\
The final version of record is available at 
\url{https://doi.org/10.1109/TAC.2018.2846684}}
\renewcommand{\IEEEarxivfooter}{Copyright (c) 2018 IEEE. Personal use is
permitted. For any other purposes, permission must be obtained from the
IEEE by emailing pubs-permissions@ieee.org.}

\makeatother

\usepackage{babel}
  \providecommand{\algorithmname}{Algorithm}
  \providecommand{\definitionname}{Definition}
  \providecommand{\propositionname}{Proposition}
  \providecommand{\remarkname}{Remark}
\providecommand{\corollaryname}{Corollary}
\providecommand{\theoremname}{Theorem}

\begin{document}
\global\long\def\interval{\mathbf{I}}

\global\long\def\realset{\ensuremath{\mathbb{R}}}

\global\long\def\rset{\ensuremath{\mathbb{R}}}

\global\long\def\complexset{\ensuremath{\mathbb{C}}}

\global\long\def\integerset{\ensuremath{\mathbb{Z}}}

\global\long\def\naturalset{\ensuremath{\mathbb{N}}}

\newcommandx\GL[2][usedefault, addprefix=\global, 1=\realset]{\ensuremath{GL\left(#2,#1\right)}}

\newcommandx\SL[2][usedefault, addprefix=\global, 1=\realset]{\ensuremath{SL\left(#2,#1\right)}}

\global\long\def\ogroup#1{O\left(#1\right)}

\global\long\def\sogroup#1{SO\left(#1\right)}

\global\long\def\SE#1{\ensuremath{SE\left(#1\right)}}

\global\long\def\SO#1{\ensuremath{SO\left(#1\right)}}

\global\long\def\SU#1{\ensuremath{SU\left(#1\right)}}

\global\long\def\sset#1{\ensuremath{S^{#1}}}

\global\long\def\ball#1#2{\ensuremath{\mathbb{B}_{#1}(#2)}}

\global\long\def\projspace#1{\ensuremath{\mathcal{P}(#1)}}

\global\long\def\sign{\ensuremath{\text{sign}}}

\global\long\def\diag#1{\text{diag}\left(#1\right)}

\global\long\def\vecop#1{\text{vec}\left(#1\right)}

\global\long\def\inprod#1#2{\left\langle #1,#2\right\rangle }

\global\long\def\abs#1{\left|#1\right|}

\global\long\def\dist{\text{dist}}

\global\long\def\norm#1{\left\Vert #1\right\Vert }

\global\long\def\grad{\ensuremath{\text{grad}}}

\global\long\def\hess{\ensuremath{\text{Hess}}}

\global\long\def\funca{f}

\global\long\def\funcb{h}

\global\long\def\funcc{g}

\global\long\def\Funca{F}

\global\long\def\axisangle#1#2{\ensuremath{R\left(#1,#2\right)}}

\global\long\def\scalar#1{#1}

\global\long\def\vector#1{#1}

\global\long\def\myvec#1{#1}

\global\long\def\veca{\myvec x}

\global\long\def\vecb{\myvec y}

\global\long\def\vecc{\myvec z}

\global\long\def\vecd{\myvec u}

\global\long\def\vecconsta{\myvec a}

\global\long\def\vecconstb{\myvec b}

\global\long\def\veccenter{\myvec c}

\global\long\def\matrixa{A}

\global\long\def\matrixb{B}

\global\long\def\matrixc{C}

\global\long\def\matrixd{D}

\global\long\def\matrixe{E}

\global\long\def\spmatrix{\matrixe}

\global\long\def\idmatrix{I}

\global\long\def\subseta{U}

\global\long\def\subsetb{U'}

\global\long\def\quatset{\ensuremath{\mathbb{H}}}
 %

\global\long\def\unitquatset{\sset 3}

\global\long\def\quat#1{\boldsymbol{#1}}

\global\long\def\imi{\hat{\imath}}

\global\long\def\imj{\hat{\jmath}}

\global\long\def\imk{\hat{k}}

\global\long\def\real#1{\mbox{Re}\left(#1\right)}

\global\long\def\imag#1{\mbox{Im}\left(#1\right)}

\global\long\def\imvec{\boldsymbol{\imath_{m}}}

\global\long\def\quatmultsymbol{\otimes}

\global\long\def\conjugsymbol{\otimes}

\global\long\def\quata{\quat q}

\global\long\def\quatb{\quat p}

\global\long\def\quatc{\quat z}

\global\long\def\dqset{\mathscr{H}}
 %

\global\long\def\udqset{\ensuremath{\dqset^{\norm 1}}}

\global\long\def\dqprimpart{\mathscr{P}}

\global\long\def\dqdualpart{\mathscr{D}}

\global\long\def\dualvector#1{\underline{\boldsymbol{#1}}}

\global\long\def\dq#1{\underline{\boldsymbol{#1}}}

\global\long\def\dual{\varepsilon}

\global\long\def\dqpta{\dq q}

\global\long\def\dqptb{\dq p}

\global\long\def\skewsymproduct#1{\ensuremath{\left\lfloor #1\right\rfloor _{\times}}}

\global\long\def\tr{\text{Tr}}

\global\long\def\diag{\text{diag}}

\global\long\def\triag{\text{triag}}

\global\long\def\cdown{\text{cdown}}

\global\long\def\vecspacea{\mathscr{V}}

\global\long\def\man#1{\boldsymbol{#1}}

\global\long\def\charta{\ensuremath{\varphi}}

\global\long\def\vecfielda{\ensuremath{\mathcal{X}}}

\global\long\def\vecfieldb{\mathcal{Y}}

\global\long\def\vecfieldc{\mathcal{Z}}

\global\long\def\vecfieldcurvea{\vecfielda}

\global\long\def\vecfieldcurveb{\vecfieldb}

\global\long\def\atlasa{\mathcal{A}}

\global\long\def\atlasb{\mathcal{B}}

\global\long\def\mana{\mathcal{N}}

\global\long\def\manb{\mathcal{R}}

\global\long\def\manc{\mathcal{\mathcal{Y}}}

\global\long\def\manfunca{f}

\global\long\def\submana{\mathcal{\mathcal{U}}}

\global\long\def\submanb{\mathcal{V}}

\global\long\def\tansubseta{V}

\global\long\def\tansubsetb{V'}

\global\long\def\tanpta{v}

\global\long\def\tanptb{u}

\global\long\def\parpta{u}

\global\long\def\manpta{\man a}

\global\long\def\manptb{\man b}

\global\long\def\manptc{\man u}

\global\long\def\manptcentral{\man c}

\global\long\def\vfieldset#1{\mathscr{X}(#1)}

\global\long\def\diffset{\mathscr{D}}

\global\long\def\identity{\text{Id}}

\global\long\def\partransp{\mbox{PT}}

\global\long\def\curve{\alpha}

\global\long\def\geod{\curve}

\global\long\def\rexp{\exp}

\global\long\def\rlog{\log}

\global\long\def\rlogb#1{\overrightarrow{#1}}

\global\long\def\arclength{\mathbb{L}}

\global\long\def\gball{\mathbb{B}}

\global\long\def\cutlocus{\mathcal{C}}

\global\long\def\tancutlocus{C}

\global\long\def\maxdefdomain{\Omega}

\global\long\def\rinj{\text{inj}}

\global\long\def\rgrad{\text{\ensuremath{\grad}}}

\global\long\def\pdf{\ensuremath{\mbox{p}}}

\global\long\def\pr{\ensuremath{\text{Pr}}}

\global\long\def\ev{\mathcal{E} }

\global\long\def\mean#1{\bar{#1}}

\global\long\def\smean{\text{\ensuremath{\mu}}}

\global\long\def\cov{P}

\global\long\def\scov{\Sigma}

\global\long\def\moment{M }

\global\long\def\smoment{\mbox{\ensuremath{\mathcal{M}}} }

\global\long\def\borelsalg{\ensuremath{\mathcal{B}}}

\global\long\def\sd{\scalar{\sigma}}

\global\long\def\ssd{s}

\global\long\def\rvset{\Phi}

\global\long\def\rva{X}

\global\long\def\rvb{Y}

\global\long\def\rvc{Z}

\global\long\def\est#1{\hat{#1}}

\global\long\def\outcome#1{\breve{#1}}

\global\long\def\normal{N}

\global\long\def\meanset{\mathbb{E}}

\global\long\def\smeanset{\mathscr{E}}

\global\long\def\state{\vector x}

\global\long\def\meas{\vector y}

\global\long\def\pnoise{\varpi}

\global\long\def\mnoise{\vartheta}

\global\long\def\pfunc{f}

\global\long\def\mfunc{h}

\global\long\def\dimstate{n_{\state}}

\global\long\def\dimmeas{n_{\meas}}

\global\long\def\dimpnoise{n_{\pnoise}}

\global\long\def\dimmnoise{n_{\mnoise}}

\global\long\def\dimstateaug{n_{a}}

\global\long\def\covpnoise{Q}

\global\long\def\covmnoise{R}

\global\long\def\vpar{v}

\global\long\def\nsp{N}

\global\long\def\weight{w}

\global\long\def\weightmatrix{W}

\global\long\def\weightm{\weight_{i}^{\left(m\right)}}

\global\long\def\weightc{\weight_{i}^{\left(c\right)}}

\global\long\def\weightcc{\weight_{i}^{\left(cc\right)}}

\global\long\def\spa{\chi}

\global\long\def\spb{\gamma}

\global\long\def\spc{\xi}

\global\long\def\spd{\zeta}

\global\long\def\sr{\sigma\text{R}}

\global\long\def\ss{\mbox{SS}}

\global\long\def\srsr{\text{SR}\sigma\text{R}}

\global\long\def\risr{\text{Ri}\sigma\text{R}}

\global\long\def\UT{\text{UT}}

\global\long\def\PaUT{\text{PaUT}}

\global\long\def\riUT{\text{RiUT}}

\global\long\def\SRUT{\text{SRUT}}

\global\long\def\riSRUT{\text{RiSRUT}}

\global\long\def\PaSRUT{\text{PaSRUT}}

\global\long\def\inov{\nu}

\global\long\def\inovquat{\quat{\nu}}

\global\long\def\kgain{G}

\global\long\def\vtq{\text{VtoQ}}

\global\long\def\RoVtoQ{\text{RoVtoQ}}

\global\long\def\GeRVtoQ{\text{GeRVtoQ}}

\global\long\def\QuVtoQ{\text{QuVtoQ}}

\global\long\def\qtv{\text{QtoV}}

\global\long\def\QtoRoV{\text{QtoRoV}}

\global\long\def\QtoGeRV{\text{QtoGeRV}}

\global\long\def\QtoQuV{\text{QtoQuV}}

\global\long\def\qwm{\text{QwMean}}

\global\long\def\fff{I}

\global\long\def\eps{eps}

\global\long\def\MT{\text{RMST}}

\global\long\def\rmsd{\text{RMSD}}

\global\long\def\samptime{\delta t}

\global\long\def\gadd{\oplus}

\global\long\def\gsub{\ominus}

\global\long\def\gsca{\odot}

\global\long\def\gprod{\otimes}

\global\long\def\biggadd{\bigoplus}

\global\long\def\biggmul{\bigotimes}

\global\long\def\gnorm#1{\norm{#1}_{\text{g}}}

\global\long\def\gid{\quat 1}

\global\long\def\ginprod#1#2{\left\langle #1,#2\right\rangle _{g}}

\global\long\def\glinset#1#2{\mathcal{L}^{g}(#1,#2)}

\global\long\def\p{\text{P}}

\global\long\def\setgrv{\Phi^{g}}

\global\long\def\gintprod#1#2{\left\langle #1,#2\right\rangle _{g}}

\global\long\def\rotmatrixa{A}

\title{Unscented Kalman Filters for Riemannian State-Space Systems}

\author{Henrique M. T. Menegaz, João Y. Ishihara, Hugo T. M. Kussaba\thanks{H. M. T. Menegaz (henriquemenegaz@unb.br) is with the Faculdade Gama at the Universidade de Brasília (UnB), Brazil.

J. Y. Ishihara (ishihara@lara.unb.br) and Hugo T. M. Kussaba (kussaba@lara.unb.br)  are with the Automation and Robotics Laboratory (LARA) at the UnB. Web-page: www.lara.unb.br.}}
\maketitle
\begin{abstract}
Unscented Kalman Filters (UKFs) have become popular in the research
community. \textit{\emph{Most UKFs work only with Euclidean }}systems,
but in many scenarios it is advantageous to consider systems with
state-variables taking values on \emph{Riemannian manifolds}. However,
we can still find some gaps in the literature's theory of UKFs for
Riemannian systems: for instance, the literature has not yet i) developed
Riemannian extensions of some fundamental concepts of the UKF theory
(e.g., extensions of $\sigma$-representation, Unscented Transformation,
Additive UKF, Augmented UKF, additive-noise system), ii) proofs of
some steps in their UKFs for Riemannian systems (e.g., proof of sigma
points parameterization by vectors, state correction equations, noise
statistics inclusion), and iii) relations between their UKFs for Riemannian
systems. In this work, we attempt to develop a theory capable of filling
these gaps. Among other results, we propose Riemannian extensions
of the main concepts in the UKF theory (including closed forms), justify
all steps of the proposed UKFs, and provide a framework able to relate
UKFs for particular manifolds among themselves and with UKFs for Euclidean
spaces. Compared with UKFs for Riemannian manifolds of the literature,
the proposed filters are more consistent, formally-principled, and
general. An example of satellite attitude tracking illustrates the
proposed theory.
\end{abstract}

\section{Introduction}

\label{sec:Introduction}

When we want to know the value of some variables of a given system–e.g.,
the position and velocity of a car, the position and attitude of a
satellite, the temperature of a boil, etc.—we can acquire data from
the system and develop a mathematical model of it. But measurements
are noisy, and models are always imperfect. Hence, to estimate the
desired variables, we often must use filters, such as Unscented Kalman
Filters (UKFs). Researchers have been applying UKFs in applications
of diverse fields: for example, in power electronic \cite{Meng2016},
aerospace \cite{Rahimi2015}, and automotive \cite{Vargas2016} systems.
These filters' success is partially explained by their good trade-off
between estimation quality and computational complexity compared with
similar techniques such as the Extended Kalman Filter (EKF) \cite{Julier2004}.

\textit{\emph{Most UKFs work only with Euclidean models (the so-called
state-space systems; cf. Section }}\ref{subsec:Additive-Unscented-Kalman}\textit{\emph{),
but sometimes modeling with }}\emph{Riemannian manifolds}\textit{\emph{
is }}better.\textit{\emph{ These manifolds can i) model more systems
(}}cf. Section \ref{subsec:Kalman-filtering-in}\textit{\emph{), ii)
provide better mathematical properties}} than Euclidean subspaces\textit{\emph{
(e.g., better metrics), and iii) be the set where measurements}} take
value from (cf. \cite{Hauberg2013,Absil2008,Pennec2006} and Section
\ref{subsec:Kalman-filtering-in}).

Although some works have introduced UKFs for Riemannian systems (e.g.,
\cite{Enayati2015,Gilitschenski2016,Lee2016,Hauberg2013}; cf. Section
\ref{subsec:Kalman-filtering-in}), we can still find some gaps in
the literature's theory for these UKFs. First, fundamental UKF concepts
still miss for Riemannian manifolds, such as $\sigma$-representation
($\sr$), Unscented Transformation (UT), Additive UKF, Augmented UKF,
additive-noise system, among others (cf. \cite{Menegaz2015}). Second,
some steps in UKFs for Riemannian manifolds are not formally justified,
such as when a UKF parameterize sigma points by vectors, or correct
the predicted state estimate, or consider noise statistics (cf. \cite{Hauberg2013,Crassidis2003,Challa2016};
see Sections \ref{subsec:Kalman-filtering-in} and \ref{subsec:Relation-with-the}).
Third, we do not know how the literature's consistent UKFs for Riemannian
manifolds relates among themselves—do they follow from a same general
Riemannian UKF?—or with UKFs for Euclidean Spaces—are these particular
cases of those?

In this work, by continuing the research of \cite{Hauberg2013}, we
aim to develop a formalized and systematized theory for UKFs on Riemannian
manifolds. Among other results, this theory introduces Riemannian
extensions of the main concepts in the UKF theory (including closed
forms), justifies all steps of the proposed UKFs, and provides a framework
able to relate UKFs for particular manifolds among themselves and
with UKFs for Euclidean spaces. 

\subsection{Kalman filtering in Riemannian manifolds}

\label{subsec:Kalman-filtering-in}

\textit{\emph{Riemannian manifolds can model many applications; far
more than Euclidean spaces. For instance, we find i)}} special orthogonal
groups, special Euclidean groups, unit spheres (including the set
of unit quaternions), and the study quadric (the set of unit dual-quaternions)
applied to many robotics applications \cite{Bullo2004,Adorno2011c,Selig2005a,Barrau2017a,Bonnabel2008},
aerospace systems \cite{Curtis2014,Wie2008,Barrau2017a,Bonnabel2008,Crassidis2003},
bio-engineering \cite{Enayati2015,Pennec2006a}, among others; ii)
positive symmetric matrices applied to applications in image recognition,
image registration, image tracking, and surgery \cite{Pennec2006a};
iii) Grassmann and Stiefel manifolds applied to information theory
\cite{Pitaval2017}, machine learning \cite{Harandi2018}, visual
recognition \cite{Harandi2018,Hajati2017}, communication systems
\cite{Seddik2017}, and geology \cite{Chepushtanova2017}; and iv)
other Riemannian manifolds applied to quantum systems \cite{AT:12},
and special and general relativity \cite{Godinho2014}.

Some works in the literature have proposed KFs for \textit{particular}
Riemannian systems: the works \cite{Condomines2013,Condomines2014,Kim2007c,Gilitschenski2016}
and \cite{Crassidis2003} (among others) introduced EKFs and UKFs
for unit quaternions; and \cite{Lee2016,DeRuiter2014} and \cite{Markovic2016}
EKFs for special orthogonal groups. Other works have proposed KFs
for \textit{classes} of Riemannian systems: the works \cite{Barczyk2011},
\cite{Martin2009} and \cite{Bonnabel2007} introduced EKFs for Lie
groups; and \cite{Hauberg2013} a UKF for geodesically-complete Riemannian
manifolds.geodesically-complete

Developing UKFs for Riemannian manifolds is difficult because, in
general, Riemannian manifolds lack some mathematical tools used in
most UKFs, such as multiplication and addition (cf. UKFs in \cite{Sarkka2013,Menegaz2016,Menegaz2015}).
An alternative is to use properties of an embedding Euclidean space
and afterwards perform operations to return to the working manifold.
For instance, an application on $\sset 3$ can use derivatives, sums,
multiplications, metrics of $\realset^{4}$ and afterwards perform
a normalization. Many works take this embedding approach \cite{Vartiainen2014,Teixeira2009,Challa2016}.

However, this approach may i) lose the physical identification (e.g.,
an addition of unit quaternions yields a non unit quaternion, which
does not represent a rotation anymore), or ii) disregard the global
properties of the manifold leading to instability. To retain the estimates
on the working manifolds, literature UKFs use intrinsic manifold properties
(cf. \cite{Barczyk2011,Martin2009,Bonnabel2007})—meaning we do not
use properties of embedding Euclidean spaces. 

In this work, we take this intrinsic approach; we combine the UKF
theory we developed in \cite{Menegaz2015} with the statistics for
Riemannian manifolds of \cite{Pennec1996} and some results of \cite{Hauberg2013}
to develop a theory of UKFs for any geodesically-complete Riemannian
manifolds. 

\section{Riemannian manifolds}

\label{sec:Riemannian-manifolds}

In this section, we provide a\textit{ }\textit{\emph{general description}}
of the concepts from Riemannian Geometry used in this work and in
Appendix \ref{appendix:Riemannian-manifolds} their \textit{\emph{formal
definitions}}. This exposition is mainly based on \cite{DoCarmo1992},
and partially on \cite{Pennec2006} and \cite{Absil2008}.

A \emph{differentiable manifold} (Definition \ref{def:differentiable-manifold})
$\mana$ (or $\mana^{n}$) can be viewed as a set whose subsets are
identified through \emph{charts} (injective mappings) with subsets
of the $\realset^{n}$. For every point $\manpta$ on a differentiable
manifold $\mana$, we can define the vector space of tangent vectors
at $\manpta$ called \textit{tangent space} and denoted by $T_{\manpta}\mana$
(Definition \ref{def:tangent-space}).

A \emph{Riemannian manifold} $\mana$ (Definition \ref{Definition:Riemannian-metric})
is a differentiable manifold endowed with a \emph{Riemannian metric}
(Definition \ref{Definition:Riemannian-metric}) $\inprod{\,}{\,}$
or $g$. For $\manpta\in\mana$ and $\tanpta\in T_{\manpta}\mana$,
the \emph{norm} of $\tanpta$ associated to $\manpta$ is defined
by $\norm{\tanpta}_{\manpta}:=\inprod{\tanpta}{\tanpta}_{\manpta}^{1/2}$
\cite{Pennec2006}. 

For two points $\manpta$ and $\manptb$ in $\mana$ connected by
a curve $\geod:\interval\rightarrow\mana$, the \emph{distance} between
$\manpta$ and $\manptb$ is defined by, for $[a,b]\subset\interval$,
\[
\dist\left(\manpta,\manptb\right):=\min_{\geod}\arclength_{a}^{b}\left(\geod\right);\quad\geod(a)=\manpta,\,\geod(b)=\manptb,
\]
where $\arclength_{a}^{b}(\curve)$ is the arc length (Definition
\ref{def:arc-length}) of $\curve$ in the interval $\left[a,b\right]$.
A \emph{geodesic} \emph{ball} of center $\man a$ and radius $r$
is the set defined as 
\[
\gball(\man a,r):=\{\man x\in\mana:\dist(\man x,\man a)<r\}.
\]

Given a tangent vector $\tanpta_{0}\in T_{\curve(t_{0})}\mana$, $t_{0}\in\interval$,
there exists only one parallel vector field $\vecfieldcurvea$ (Definition
\ref{def:vector-field}) along $\curve$, such that $\vecfieldcurvea(t_{0})=\tanpta_{0}$;
$\vecfieldcurvea(t)$ is called the \emph{parallel transport} of $\vecfieldcurvea(t_{0})$
along $\curve$.

A curve $\geod:\interval\rightarrow\mana$ is called a \emph{geodesic
at} $t_{0}\in\interval$ if 
\[
\frac{D}{dt}\Big(\geod'(t)\Big)=0
\]
at $t_{0}$, where $D/dt(\geod'(t))$ is the covariant derivative
of $\geod'(t)$ (Theorem \ref{thm:covariant-derivative-1}); if $\geod$
is a geodesic at $t$, for all $t\in\interval$, we say $\geod$ is
a \emph{geodesic} \cite{DoCarmo1992}. If a curve minimizes the \emph{arc
length} between two points of the manifold, then this curve is a geodesic,
but the converse is \emph{only} valid \emph{locally}. If the definition
domain of all geodesics of $\mana$ can be extended to $\realset$,
then $\mana$ is said to be \emph{geodesically-complete}. There exists
at least one geodesic connecting every two points of a geodesically-complete
manifold.

Given a point $\manpta\in\mana$, the \emph{exponential mapping} (Definition
\ref{def:exponential-map}), denoted by $\rexp_{\manpta}$, associates
a vector of $T_{\manpta}\mana$ to a point of $\mana$. Geometrically,
$\rexp_{\manpta}(\tanpta)$ is a point of $\mana$ obtained by going
out the length equal to $\norm{\tanpta}$, starting from $\manpta$,
along a geodesic which passes through $\manpta$ with velocity equal
to $v/\norm{\tanpta}$. 

Assuming a geodesically-complete manifold, it is possible to follow
the geodesic $\rexp_{\manpta}(t\tanpta)$ from $t=0$ to $t\rightarrow\infty$.
It may happen, however, that from a particular value $t_{\tanpta}$
to $t\rightarrow\infty$, the geodesics $\rexp_{\manpta}(t\tanpta)$
\emph{do not minimize the arc length} between $\manpta$ and $\rexp_{\manpta}(t\tanpta)$.
In this case, the subset $\{\rexp_{\manpta}(t_{v}\tanpta):\tanpta\in T_{\manpta}\mana\}\subset\mana$
is called the \emph{cut locus} $\cutlocus(\manpta)$ and the inverse
image $\tancutlocus(\manpta):=\rexp_{\manpta}^{-1}[\cutlocus(\manpta)]$
the \emph{tangential cut locus} \cite{Pennec2006}. The \emph{injectivity
radius} of $\mana$ is defined as $\rinj(\mana):=\inf_{\man p\in\mana}\dist(\man p,\cutlocus(\man p)).$ 

For every $\manpta\in\mana$, we can reduce the domain of $\rexp_{\manpta}$
to some subsets such that $\rexp_{\manpta}$ is a diffeomorphism.
The maximal of these subsets is called the \emph{maximal definition
domain} $\maxdefdomain(\manpta)\subset T_{\manpta}\mana$; this set
is \emph{bounded} by $\tancutlocus(\manpta)$ \cite{Pennec2006}.
The inverse mapping of $\rexp_{\manpta}$ is the (Riemannian) \emph{logarithm
mapping} (Definition \ref{def:exponential-map})\textit{\emph{ and}}
we denote it by either \emph{$\rlog_{\manpta}\manptb$} or $\overrightarrow{\manptb\manpta}$.

\section{Intrinsic Statistics on Riemannian Manifolds}

\label{sec:Intrinsic-Statistics-on}

UKFs are based on information of \emph{moments} of \emph{random vectors}
and of \emph{sample moments} of \emph{weighted sets}. To define UKFs
on Riemannian manifolds, we need extensions of these concepts.

\subsection{Statistics of random points}

\label{subsec:Statistics-of-random}

Riemannian extensions of random vectors are called \emph{(Riemannian)
random points} \cite{Pennec2006}; the set of all random points taking
values on a Riemannian manifold $\mana$ is denoted by $\mbox{\ensuremath{\man{\rvset}}}_{\mana}$.
Given a random point $\man{\rva}\in\mbox{\ensuremath{\man{\rvset}}}_{\mana}$,
its \emph{probability density function (pdf)} is denoted by $\man{\pdf}_{\man{\rva}}$,
and for a real-valued function $\Funca:\mana\rightarrow\realset$
the\emph{ expected value of $\Funca$ relative to $\man{\rva}$ is
defined by
\begin{equation}
\man{\ev}_{\man{\rva}}\left\{ F(\man{\rva})\right\} :=\int_{\mana}\Funca(\man{\manptb})\man{\pdf}_{\man{\rva}}(\man{\manptb})d\mana(\man{\manptb}).\label{eq:expectation-of-real-valued-function}
\end{equation}
}For functions taking values on manifolds, we cannot define the expected
value as in (\ref{eq:expectation-of-real-valued-function}); thus,
we define mean points following the Karcher expectation: they are
the \emph{local} minima of variances \cite{Pennec2006}. 

Given a point $\manptcentral\in\mana$, the \emph{variance} $\sd_{\man{\rva}}^{2}(\manptcentral)$
is defined by $\sd_{\man{\rva}}^{2}(\manptcentral):=\man{\ev}_{\man{\rva}}\{\dist^{2}(\manptcentral,\man{\rva})\}.$
If $\sd_{\man{\rva}}^{2}(\man c)$ is finite for every point $\manptcentral\in\mana$,
then a point $\mean{\man{\rva}}\in\mana$ is an \emph{expected point}
or \emph{mean} of $\man{\rva}$ if 
\begin{equation}
\mean{\man{\rva}}=\arg\underset{\man c\in\mana}{\min}\sd_{\man{\rva}}^{2}(\manptcentral).\label{eq:Riemannian-mean-definition}
\end{equation}
The set of all means of $\man{\rva}$ is denoted by $\meanset(\man{\rva})$.
A random point can have more than one mean\footnote{For a discussion about the existence and uniqueness of this expectation,
cf. Section 4.2 of \cite{Pennec2006}.}.

Let $\man{\rva}\in\mbox{\ensuremath{\man{\rvset}}}_{\mana}$ be a
random point with a mean $\mean{\man{\rva}}\in\meanset(\man{\rva})$,
and consider a point $\manpta\in\mana$. If $\mean{\man{\rva}}\in\maxdefdomain(\manpta)$,
then the \emph{$j$th (central) moment of $\man{\rva}$ with respect
to $\mean{\man{\rva}}$ at $\manpta$} is defined by, for even $j$,
\begin{equation}
\man{\moment}_{\man{\rva},\mean{\man{\rva}}}^{\manpta,j}:=\man{\ev}_{\man{\rva}}\Big\{\Big[\Big(\overrightarrow{\manpta\man{\rva}}-\rlogb{\manpta\mean{\man{\rva}}}\Big)(\diamond)^{T}\Big]^{\otimes\frac{j}{2}}\Big\};\label{eq:Riemannian-central-moment}
\end{equation}
and for odd $j$,
\[
\man{\moment}_{\man{\rva},\mean{\man{\rva}}}^{\manpta,j}:=\man{\ev}_{\man{\rva}}\Big\{\Big[\Big(\overrightarrow{\manpta\man{\rva}}-\rlogb{\manpta\mean{\man{\rva}}}\Big)(\diamond)^{T}\Big]^{\otimes\frac{j-1}{2}}\otimes\Big(\overrightarrow{\manpta\man{\rva}}-\rlogb{\manpta\mean{\man{\rva}}}\Big)\Big\}.
\]
We define joint pdf {[}denoted by $\man{\pdf}_{\man{\rva}\man Y}(\man x,\man y)${]},
joint expected moment ($\man{\ev}_{\man{\rva}\man Y}\left\{ \manfunca(\man x,\man y)\right\} $)
and cross-covariance ($\man P_{\man{\rva}\man Y,(\mean{\man{\rva}},\bar{\man Y})}^{\manpta\man b}$)
of two random points $\man{\rva}$ and $\man Y$ similarly (cf. \cite{Menegaz2016}).
The notation $\man{\rva}\sim(\mean{\man{\rva}},\man{\moment}_{\man{\rva},\mean{\man{\rva}}}^{\manpta,2},...,\man{\moment}_{\man{\rva},\mean{\man{\rva}}}^{\manpta,l})_{\mana}$
stands for a Riemannian random point $\man{\rva}\in\man{\Phi}_{\mana}$
with mean $\mean{\man{\rva}}\in\meanset(\man{\rva})$ and moments
$\man{\moment}_{\man{\rva},\mean{\man{\rva}}}^{\manpta,2},...,\man{\moment}_{\man{\rva},\mean{\man{\rva}}}^{\manpta,l}$.
The second moment ($j=2$) is called \emph{covariance} and denoted
by $\man P_{\man{\rva}\man{\rva},\mean{\man{\rva}}}^{\manpta}:=\man{\moment}_{\man{\rva},\mean{\man{\rva}}}^{\manpta,j}$.
If $\meanset(\man{\rva})=\{\mean{\man{\rva}}\},$ we can write $\man{\moment}_{\man{\rva}}^{\manpta,j}:=\man{\moment}_{\man{\rva},\mean{\man{\rva}}}^{\manpta,j}$
and $\man P_{\man{\rva}\man{\rva}}^{\manpta}:=\man P_{\man{\rva}\man{\rva},\mean{\man{\rva}}}^{\manpta}$,
or even $\man{\moment}_{\man{\rva}}^{j}:=\man{\moment}_{\man{\rva},\mean{\man{\rva}}}^{\mean{\man{\rva}},j}$
and $\man P_{\man{\rva}\man{\rva}}:=\man P_{\man{\rva}\man{\rva},\mean{\man{\rva}}}^{\mean{\man{\rva}}}$.

We represent statistics of \emph{Euclidean} manifolds without bold
notation. For $\rva\in\rvset_{\realset^{n}}$, $X$ is symmetric if
{\small{}$\pdf_{X}(\bar{X}+x)=\pdf_{X}(\bar{X}-x)$} for every $x\in\mathbb{R}^{n}$.
If $\rva$ has a mean, then 
\[
\mean{\rva}=\arg\underset{\veccenter\in\realset^{n}}{\min}\sd_{\rva}^{2}(\veccenter)=\ev_{\rva}\{X\};
\]
and, for $j$ even,
\[
\moment_{X}^{j}:=\ev_{\rva}\{[(\overrightarrow{\mean{\rva}\rva}-\rlogb{\mean{\rva}\mean{\rva}})(\diamond)^{T}]^{\otimes\frac{j}{2}}\}=\ev_{\rva}\{[(\rva-\mean{\rva})(\diamond)^{T}]^{\otimes\frac{j}{2}}\}
\]
(similarly for $j$ odd and for sample cross-covariances).

\subsection{Statistics of weighted sets}

\label{sec:Statistics-of-weighted}

For a Riemannian manifold $\mana$ and the natural numbers $l\geq2$
and $\nsp\geq1$, consider the weighted set
\begin{multline*}
\man{\chi}:=\Big\{\man{\spa}_{i},\weight_{i}^{m},\weight_{i}^{c,j},\weight_{i}^{cc,j}:\man{\chi}_{i}\in\mana;\\
j=1,...,l;\,\weight_{i}^{m},\weight_{i}^{c,j},\weight_{i}^{cc,j}\in\realset\Big\}_{i=1}^{N}.
\end{multline*}
The weights $w_{i}^{m}$ are associated (below) with the definition
of sample\emph{ }\textit{\emph{mean}}, $w_{i}^{c,j}$ with the $j$th
sample\emph{ }\textit{\emph{moment}}, and $w_{i}^{cc}$ with the $j$th
sample \textit{\emph{cross-moment}} of $\man{\chi}$.

The\emph{ sample variance }of $\man{\chi}$ with respect to a point
$\man c\in\mana$ is defined by $\ssd_{\mathbb{\man{\chi}}}^{2}(\manptcentral):=\sum_{i=1}^{N}w_{i}^{m}\dist^{2}(\manptcentral,\man{\chi}_{i}).$
If the variance $\ssd_{\mathbb{\man{\chi}}}^{2}(\manptcentral)$ is
finite for every point $\manptcentral\in\mana$, then a \emph{sample}
\emph{expected point} or\emph{ sample mean} of $\man{\chi}$ is defined
by
\begin{equation}
\man{\mu}_{\man{\chi}}:=\arg\underset{\man c\in\mana}{\min}\ssd_{\mathbb{\man{\chi}}}^{2}(\manptcentral).\label{eq:Riemannian-sample-mean-definition}
\end{equation}
The set of all sample means of $\man{\chi}$ is represented by $\smeanset(\man{\chi})$.
An weighted set in the form of $\man{\chi}$ can have more than one
sample mean.

For a point $\manpta\in\mana$, if $\man{\mu}_{\man{\chi}},\man{\chi}_{1},\man{\chi}_{2},...,\man{\chi}_{N}\in\mana-\cutlocus(\manpta)$,
then the\emph{ $j$th sample moment of $\man{\chi}$ with respect
to $\mean{\man{\rva}}$ at $\manpta$} is defined by, for $j$ even,
\begin{equation}
\man{\smoment}_{\man{\chi},\man{\mu}_{\man{\chi}}}^{\manpta,j}:=\sum_{i=1}^{N}\weight_{i}^{c,j}\Big[\big(\rlogb{\manpta\man{\chi}_{i}}-\rlogb{\manpta\man{\mu}_{\man{\chi}}}\big)(\diamond)^{T}\Big]^{\otimes\frac{j}{2}};\label{eq:Riemannian-sample-moment-definition}
\end{equation}
and for $j$ odd,
\[
\man{\smoment}_{\man{\chi},\man{\mu}_{\man{\chi}}}^{\manpta,j}:=\sum_{i=1}^{N}\weight_{i}^{c,j}\Big[\big(\rlogb{\manpta\man{\chi}_{i}}-\rlogb{\manpta\man{\mu}_{\man{\chi}}}\big)(\diamond)^{T}\Big]^{\otimes\frac{j-1}{2}}\otimes\big(\rlogb{\manpta\man{\chi}_{i}}-\rlogb{\manpta\man{\mu}_{\man{\chi}}}\big).
\]
The sample moment ($j=2$) is called \emph{sample covariance} and
denoted by $\man{\scov}_{\man{\chi}\man{\chi},\man{\mu}_{\man{\chi}}}^{\manpta}:=\man{\smoment}_{\man{\chi},\man{\mu}_{\man{\chi}}}^{\manpta,2}$.
If $\smeanset\left(\man{\chi}\right)=\{\man{\mu}_{\man{\chi}}\},$
we can write $\man{\smoment}_{\man{\chi}}^{\manpta,j}:=\man{\smoment}_{\man{\chi},\man{\mu}_{\man{\chi}}}^{\manpta,j}$
and $\man{\scov}_{\man{\chi}\man{\chi}}^{\manpta}:=\man{\scov}_{\man{\chi}\man{\chi},\man{\mu}_{\man{\chi}}}^{\manpta}$;
or even, $\man{\smoment}_{\man{\chi}}^{j}:=\man{\smoment}_{\man{\chi},\man{\mu}_{\man{\chi}}}^{\man{\mu}_{\man{\chi}},j}$and
$\man{\scov}_{\man{\chi}\man{\chi}}:=\man{\scov}_{\man{\chi}\man{\chi},\man{\mu}_{\man{\chi}}}^{\man{\mu}_{\man{\chi}}}$. 

In addition, for i) the Riemannian manifold $\manb$, ii) a function
$f:\mana\rightarrow\manb$, iii) the weighted set
\[
\mathbf{\man{\spb}}:=\Big\{\man{\spb}_{i},\weight_{i}^{m},\weight_{i}^{c,j},\weight_{i}^{cc,j}:\man{\spb}_{i}=f(\man{\spa}_{i});\,j=1,...,l\Big\}_{i=1}^{\nsp},
\]
with a mean $\man{\mu}_{\man{\spb}}$, and iv) the point $\manptb\in\manb$.
If $\man{\mu}_{\man{\spb}},\man{\spb}_{1},\man{\spb}_{2},...,\man{\spb}_{N}\in\manb-\cutlocus(\manptb)$,
then the \emph{$j$th} \emph{cross-moment of $\man{\chi}$ and $\man{\spb}$
with respect to $(\man{\mu}_{\man{\spa}},\man{\mu}_{\man{\spb}})$
at $(\manpta,\manptb)$} is defined by, for $j$ even,
\[
\man{\smoment}_{\man{\spa}\man{\spb},\man{\mu}_{\man{\spa}}\man{\mu}_{\man{\spb}}}^{j,\manpta\manptb}:=\sum_{i=1}^{N}\weight_{i}^{cc,j}\Big[\big(\rlogb{\manpta\man{\spa}_{i}}-\rlogb{\manpta\man{\mu}_{\man{\spa}}}\big)\big(\rlogb{\manpta\man{\spb}_{i}}-\rlogb{\manpta\man{\mu}_{\man{\spb}}}\big)^{T}\big)\Big]^{\otimes\frac{j}{2}};
\]
and for $j$ odd,
\begin{multline*}
\man{\smoment}_{\man{\spa}\man{\spb},\man{\mu}_{\man{\spa}}\man{\mu}_{\man{\spb}}}^{j,\manpta\manptb}:=\sum_{i=1}^{N}\weight_{i}^{cc,j}\Big[\big(\rlogb{\manpta\man{\spa}_{i}}-\rlogb{\manpta\man{\mu}_{\man{\spa}}}\big)\\
\times\big(\rlogb{\manpta\man{\spb}_{i}}-\rlogb{\manpta\man{\mu}_{\man{\spb}}}\big)^{T}\big)\Big]^{\otimes\frac{j-1}{2}}\otimes\big(\rlogb{\manpta\man{\spa}_{i}}-\rlogb{\manpta\man{\mu}_{\man{\spa}}}\big).
\end{multline*}
The second sample cross-moment ($j=2$) is called \emph{sample cross-covariance}
and denoted by $\man{\scov}_{\man{\chi}\man{\spb},\man{\mu}_{\man{\spa}}\man{\mu}_{\man{\spb}}}^{\manpta\manptb}:=\man{\smoment}_{\man{\spa}\man{\spb},\man{\mu}_{\man{\spa}}\man{\mu}_{\man{\spb}}}^{j,\manpta\manptb}$.
If $\smeanset\left(\mathbf{\man{\spa}}\right)=\{\man{\mu}_{\man{\spa}}\}$
and $\smeanset\left(\mathbf{\man{\spb}}\right)=\{\man{\mu}_{\man{\spb}}\},$
we can write $\man{\smoment}_{\man{\spa}\man{\spb}}^{j,\manpta\manptb}:=\man{\smoment}_{\man{\spa}\man{\spb},\man{\mu}_{\man{\spa}}\man{\mu}_{\man{\spb}}}^{j,\manpta\manptb}$
and $\man{\scov}_{\man{\chi}\man{\spb}}^{\manpta\manptb}:=\man{\scov}_{\man{\chi}\man{\spb},\man{\mu}_{\man{\spa}}\man{\mu}_{\man{\spb}}}^{\manpta\manptb}$;
or even, if $\man{\smoment}_{\man{\spa}\man{\spb}}^{j}:=\man{\smoment}_{\man{\spa}\man{\spb}}^{j,\man{\mu}_{\man{\spa}}\man{\mu}_{\man{\spb}}}$and
$\man{\scov}_{\man{\chi}\man{\spb}}:=\man{\scov}_{\man{\chi}\man{\spb}}^{\man{\mu}_{\man{\spa}}\man{\mu}_{\man{\spb}}}$. 

We represent \emph{Euclidean} sets sample statistics without bold
notation. For a set $\chi$ with points $\chi_{i}\in\realset^{n}$,
we have
\[
\smean_{\spa}=\arg\underset{\veccenter\in\realset^{n}}{\min}\ssd_{\mathbb{\chi}}^{2}(\veccenter)=\sum_{i=1}^{n}\weight_{i}^{m}\chi_{i};
\]
and, for $j$ even, $\smoment_{\chi}^{j}=\sum_{i=1}^{N}w_{i}^{c,j}[(\chi_{i}-\mu_{\chi})(\diamond)^{T}]^{\otimes\frac{j}{2}}$
(similarly for $j$ odd and for sample cross-moments).

\section{Unscented Kalman Filters}

\label{subsec:Additive-Unscented-Kalman}

There are two main concepts required to define UKFs, namely: $\sr$s
and $\UT$s \cite{Menegaz2015}. Broadly, i) a $\sr$ is a set of
weighted points (the sigma points) approximating a random vector,
and ii) a UT is a function mapping two functionally related random
vectors to two sets that approximate their joint pdf.

For the natural numbers $l\geq2$ and $\nsp\geq1$, consider i) a
function $f:\realset^{n}\rightarrow\realset^{\eta}$; ii) the random
vectors $\rva\sim(\mean{\rva},\moment_{\rva}^{2},...,\moment_{\rva}^{l})_{\realset^{n}}$
and $\rvb:=\funca(\rva)\sim(\mean{\rvb},\moment_{\rvb}^{2},...,\moment_{\rvb}^{l})_{\realset^{\eta}}$;
and iii) the sets\footnote{Compared with \cite{Menegaz2015}, here we consider \emph{simpler}
sets. With this consideration, we have a clearer text and \emph{do
not} \emph{lose generality} for the results relative to the UKFs.} 
\begin{multline*}
\chi:=\Big\{\spa_{i},\weight_{i}^{m},\weight_{i}^{c,j},\weight_{i}^{cc,j}:\chi_{i}\in\realset^{n};\\
j=1,...,l;\,\weight_{i}^{m},\weight_{i}^{c,j},\weight_{i}^{cc,j}\in\realset\Big\}_{i=1}^{N};\text{ and}
\end{multline*}
\[
\spb:=\Big\{\spb_{i},\weight_{i}^{m},\weight_{i}^{c,j},\weight_{i}^{cc,j}:\spb_{i}=f(\spa_{i});\,j=1,...,l\Big\}_{i=1}^{\nsp}.
\]
\begin{defn}[$\sigma$R. Definition 1 of \cite{Menegaz2015}]
\label{def:sigma_representacao_ext1} The set $\spa$ is an $l$\emph{th
order} $\nsp$ \emph{points} $\sigma$R\emph{($l$th$\nsp\sigma$R)
of} $\rva$ if, for every $j=1,\ldots,l$:
\begin{align}
\weight_{i}^{m} & \neq0,\,\weight_{i}^{c,j}\neq0,\weight_{i}^{cc,j}\neq0,\quad i=1,\ldots,\nsp;\label{eq:sigma-rep-def-weights-condition}\\
\smean_{\spa} & =\mean{\rva};\label{eq:sigma-rep-def-mean-condition}\\
\smoment_{\spa}^{j} & =\moment_{\rva}^{j}.\label{eq:sigma-rep-def-moments-condition}
\end{align}
\end{defn}
\begin{defn}[UT. Definition 2 of \cite{Menegaz2015}]
\label{def:unscented-transform} If $\smean_{\spa}=\mean{\rva}$
and $\smoment_{\spa}^{j}=\moment_{\rva}^{j}$ for every $j=2,\ldots,l$;
then the $l$\emph{th order UT ($l$UT)} is defined by 
\begin{multline*}
l\mbox{UT}:\big(f,\mean{\rva},\moment_{\rva}^{2},...,\moment_{\rva}^{l}\big)\mapsto\\
\big(\smean_{\spb},\smoment_{\spb}^{2},...,\smoment_{\spb}^{l},\smoment_{\spa\spb}^{2},...,\smoment_{\spa\spb}^{l}\big).
\end{multline*}
$\spa$ is called the \emph{independent set} of an $l\mbox{UT}$,
and $\spb$ its \emph{dependent} set.
\end{defn}
Every $l$th$\nsp\sigma$R is an independent set of an $l\mbox{UT}$.
When calling an $l$th$N\sigma$R of $X$ or an $l\UT$, the reference
to the $l$th order can be omitted if $l=2$. Also, the reference
to $N$ point and/or to $X$ can be omitted in case they are obvious
from the context or irrelevant to a discussion.

We can apply UTs in KF prediction-correction frameworks to form UKFs.
UKFs estimate the state of systems described either in the additive
form 
\begin{equation}
\state_{k}=\pfunc_{k}\left(\state_{k-1}\right)+\pnoise_{k},\,\meas_{k}=\mfunc_{k}\left(\state_{k}\right)+\mnoise_{k};\label{eq:additive-system}
\end{equation}
or, more generally, in the form 
\begin{equation}
\state_{k}=\pfunc_{k}\left(\state_{k-1},\pnoise_{k}\right),\,\meas_{k}=\mfunc_{k}\left(\state_{k},\mnoise_{k}\right),\label{eq:general-system}
\end{equation}
where $k$ is the time step; $\state_{k}$ $\in\Phi^{\dimstate}$
is the internal state; $\meas_{k}\in\Phi^{\dimmeas}$ is the measured
output; and $\pnoise_{k}\in\Phi^{\dimpnoise}$ and $\mnoise_{k}\in\Phi^{\dimmnoise}$
are the process and measurement noises respectively; the noise terms
$\pnoise_{k}$ and $\mnoise_{k}$ are assumed to be uncorrelated.

In \cite{Menegaz2015}, we developed consistent UKFs for these systems:
the the\textit{ }\emph{Additive UKF}\textit{\emph{ (}}\emph{AdUKF,
}Algorithm 6 of \cite{Menegaz2016}; see also \cite{Menegaz2015})
for (\ref{eq:additive-system}); and the \textit{Augmented UKF} (\textit{AuUKF},
Algorithm 7 of \cite{Menegaz2016}; see also \cite{Menegaz2015})
for (\ref{eq:general-system}). But how could we develop similar UKFs
when $\state_{k}$, $\meas_{k}$, $\pnoise_{k}$ and $\mnoise_{k}$
are Riemannian random points? In the next section, we begin a theory
towards this goal. 

\section{Riemannian $\sigma$-representations}

\label{sec:Riemannian-sigma-representations}

In this section, first, we define \emph{Riemannian $\sigma$-representation}s
($\risr$). They extend $\sr$s to Riemannian manifolds: $\sr$s approximate
random vectors, and $\risr$s approximate Riemannian random points.
Then, we show a way of extending closed forms of $\sr$s to $\risr$s.
Afterwards, we introduce results relative to the minimum number of
sigma points of an $\risr$. At last, we introduce some particular
forms of $\risr$s.

For now on, we make the following assumptions—we explain their implications
in Section \ref{subsec:Riemannian-Unscented-Filters}—:
\begin{enumerate}
\item \label{enu:assumption1}all Riemannian manifolds are \emph{geodesically-complete};
\item \label{enu:assumption2}all Riemannian exponential mappings are defined
with their domain allowing them to \emph{realize diffeomorphisms;}
\item \label{enu:assumption6}every set of weighted points belonging to
a Riemannian manifold admits \emph{one, and only one, Riemannian sample
mean.}
\end{enumerate}
For the point $\manpta\in\mana$ and the natural numbers $l\geq2$
and $\nsp\geq1$, consider i) a random point $\man{\rva}\sim(\mean{\man{\rva}},\man{\moment}_{\man{\rva},\mean{\man{\rva}}}^{\manpta,2},...,\man{\moment}_{\man{\rva},\mean{\man{\rva}}}^{\manpta,l})_{\mana^{n}}$
and ii) a weighted set $\man{\spa}:=\{\man{\spa}_{i},\weight_{i}^{m},\weight_{i}^{c,j},\weight_{i}^{cc,j}|\man{\spa}_{i}\in\mana\}{}_{i=1}^{\nsp}$
with sample mean $\man{\smean}_{\man{\spa}}$ and sample moments $\man{\smoment}_{\man{\spa}}^{j}$,
$j=2$, ..., $l$. 

\begin{defn}[\emph{Ri$\sr$. Definition 9.1 of \cite{Menegaz2016}}]
\label{def:Consider-a-Riemannian}The set $\man{\spa}$ is a \emph{Riemannian}
$l$\emph{th order} $\nsp$ \emph{points} $\sigma$-\emph{representation
(Ri$l$th$\nsp\sigma$R) of} $\man{\rva}$ if, for every $j=1,\ldots,l$:
\begin{align}
\weight_{i}^{m} & \neq0,\,\weight_{i}^{c,j}\neq0,\weight_{i}^{cc,j}\neq0,\quad i=1,\ldots,\nsp;\label{eq:Riemannian-sigma-rep-definition-condition}\\
\man{\smean}_{\man{\spa}} & =\mean{\man{\rva}};\label{eq:Riemannian-sigma-rep-definition-mean-condition}\\
\man{\smoment}_{\man{\spa}}^{j} & =\man{\moment}_{\man{\rva}}^{j},\quad j=2,3,\ldots,l;\label{eq:Riemannian-sigma-rep-definition-moments-condition}
\end{align}

Moreover, assume $\man{\spa}$ is an Ri$l$th$\nsp\sigma$R of $\rva$,
then: 
\begin{itemize}[labelsep=0.1cm,leftmargin=0.30cm]
\item $\spa$ is \emph{normalized} if, for every $j=1,2,\ldots,l$:
\[
\sum_{i=1}^{\nsp}\weight_{i}^{m}=\sum_{i=1}^{\nsp}\weight_{i}^{c,j}=\sum_{i=1}^{\nsp}\weight_{i}^{cc,j}=1.
\]
\item $\spa$ is \emph{homogeneous} if, for every $j=1,2,\ldots,l$, the
following equations are satisfied: for $N$ odd and every $i=1,...,\nsp-1$:
\begin{equation}
\weight_{1}^{m}=\weight_{i}^{m},\,\weight_{1}^{c,j}=\weight_{i}^{c,j},\,\weight_{1}^{cc,j}=\weight_{i}^{cc,j};\label{eq:Riemannian-sigma-rep-def-homogeneous-odd}
\end{equation}
or, for $N$ even and every $i=1,...,\nsp$:
\begin{equation}
\weight_{1}^{m}=\weight_{i}^{m},\,\weight_{1}^{c,j}=\weight_{i}^{c,j},\,\weight_{1}^{cc,j}=\weight_{i}^{cc,j}.\label{eq:Riemannian-sigma-rep-def-homogeneous-even}
\end{equation}
\item $\man{\spa}$ is \emph{symmetric (with respect to $\man{\spa}_{\nsp}$,
without loss of generality)} if
\begin{multline}
\rlogb{\man{\smean}_{\man{\spa}}\man{\spa}_{i}}-\rlogb{\man{\smean}_{\man{\spa}}\man{\spa}_{\nsp}}=-\Big(\rlogb{\man{\smean}_{\man{\spa}}\man{\spa}_{i+\text{int}\left(\frac{\nsp}{2}\right)}}-\rlogb{\man{\smean}_{\man{\spa}}\man{\spa}_{\nsp}}\Big),\\
\weight_{i}^{m}=\weight_{i+\text{int}\left(\frac{N}{2}\right)}^{m},\,\weight_{i}^{c,j}=\weight_{i+\text{int}\left(\frac{N}{2}\right)}^{c,j},\,\weight_{i}^{cc,j}=\weight_{i+\text{int}\left(\frac{N}{2}\right)}^{cc,j},\label{eq:Riemannian-sigma-rep-def-symmetric}
\end{multline}
for every $j=1,2,\ldots,l$ and $i=1,...,\text{int}(\nsp/2)$, where
$\text{int}(\nsp/2)$ stands for greatest integer less than or equal
to $\nsp/2$. 
\end{itemize}
\end{defn}
When calling an Ri$l$th$\nsp\sigma$R of $\man{\rva}$, the reference
to the $l$th order can be omitted if $l=2$. Also, the reference
to $\nsp$ points or to $\man{\rva}$ can be omitted if they are obvious
from the context or irrelevant to a discussion. 

Ri$l$th$\nsp\sigma$Rs are generalizations of $l$th$\nsp\sigma$Rs;
every $l$th$\nsp\sigma$R with an Ri$l$th$\nsp\sigma$R, and every
Ri$l$th$\nsp\sigma$R with Euclidean points is an $l$th$\nsp\sigma$R.
This follows directly from the last paragraph of Sections \ref{subsec:Statistics-of-random}
and of \ref{sec:Statistics-of-weighted}.

Finding closed forms for $\risr$s may be troublesome, but the next
theorem provides a way of obtaining them from closed forms of $\sr$s—the
reader will find several closed forms of $\sr$s in \cite{Menegaz2015,Menegaz2016,Sarkka2013}.

\begin{thm}[\emph{Theorem 9.1 of \cite{Menegaz2016}}]
\label{thm:Euclidean-to-Riemannian-sigma-rep}Suppose that, for every
$i=1,\ldots,\nsp$,
\begin{enumerate}
\item $\weight_{i}^{m}>0$, 
\item $\maxdefdomain(\man{\mean{\rva}})$ is convex, and 
\item $\man{\spa}_{i}\in\gball(\man{\mean{\rva}},r)\cap\cutlocus(\man{\mean{\rva}})$
\end{enumerate}
where $0<r\le\frac{1}{2}\min\{\rinj(\mana),\pi/\sqrt{\kappa})$ and
$\kappa$ is an upper bound of the sectional curvatures of $\mathcal{N}$.
Then \textup{$\man{\spa}$} is a normalized Ri$l$th$\nsp\sigma$R
of $\man{\rva}$ if, and only if, 
\[
\spa:=\big(\rlog_{\man{\mean{\rva}}}\man{\spa}_{i},\weight_{i}^{m},\weight_{i}^{c,j},\weight_{i}^{cc,j}\big)_{i=1}^{\nsp}
\]
is a normalized $l$th$\nsp\sigma$R of the random vector 
\[
\rva\sim\big([0]_{n\times1},\man{\moment}_{\man{\rva}}^{2},\ldots,\man{\moment}_{\man{\rva}}^{l}\big)_{T_{\man{\mean{\rva}}}\mana}.
\]
Moreover, the following statements are true:

\begin{enumerate}
\item \label{enu:theo:euclidean-to-Riemannian-sigma-rep-homoegeneous}$\man{\spa}$
is homogeneous if, and only if, $\spa$ is homogeneous;
\item \label{enu:theo:euclidean-to-Riemannian-sigma-rep-symmetric}$\man{\spa}$
is symmetric if, and only if, $\spa$ is symmetric.
\end{enumerate}
\end{thm}

The proof of Theorem \ref{thm:Euclidean-to-Riemannian-sigma-rep}
is given in Appendix \ref{proof:theoremRiSRs}; for conditions to
assure the convexity of $\maxdefdomain(\man{\mean{\rva}})$, see \cite{FGR:15}
and references therein.

With this theorem, we can extend some results from $l$th$\nsp\sigma$Rs
to Ri$l$th$\nsp\sigma$Rs, such as the minimum number of sigma points
of an Ri$l$th$\nsp\sigma$R. 

\begin{cor}[Corollary 9.1 of \cite{Menegaz2016}]
\label{cor:Riemannian-sr-minimum-numbers}Let i) $\man{\spa}$ be
a normalized Ri$l$th$\nsp\sigma$R of $\man{\rva}$ with $\weight_{i}^{m}>0$
for every $i=1,\ldots,\nsp$; and ii) the rank of the covariance $\man{\cov_{\rva\rva}}$
be $r\leq n$. Then the following statements are true: 
\begin{enumerate}
\item \label{enu:Minimum number of sigma points1-2}$\nsp\geq r+1$. If
$\nsp=r+1$, then $\man{\spa}$ is called a \emph{minimum} Ri$l$th$\nsp\sigma$R
of $\man{\rva}$. 
\item \label{enu:Minimum number of sigma points1-1-1}If $\man{\spa}$ is
symmetric, then $\nsp\geq2r$. If $\man{\spa}$ is symmetric and $\nsp=2r$,
then $\man{\spa}$ is called a \emph{minimum symmetric} Ri$l$th$\nsp\sigma$R
of $\man{\rva}$. 
\end{enumerate}
Moreover, consider the set $\spa:=\big\{\rlogb{\mean{\man{\rva}}\man{\spa}_{i}},\weight_{i}^{m},\weight_{i}^{c,j},\weight_{i}^{cc,j}\big\}_{i=1}^{\nsp}$
and the random vector $\rva\sim\big([0]_{n\times1},\man{\cov}_{\man{\rva}\man{\rva}}\big)_{T_{\mean{\man{\rva}}}\mana}.$
Then the following statements are true: 
\begin{itemize}
\item If $\spa$ is a (normalized) homogeneous minimum symmetric $\sr$
of $\rva$ (\emph{HoMiSy$\sigma$R}, Corollary 3 of \cite{Menegaz2015}),
then $\man{\spa}$ is also minimum and symmetric and is called a \emph{Riemannian
(normalized) homogeneous minimum symmetric $\sigma$ -representation
of $\man{\rva}$.}
\item If $\spa$ is a \emph{Rho Minimum $\sigma$R} of $\rva$ (``it is
described in the 6th row of Table I of \cite{Menegaz2015} and refereed
there as the “Minimum set of {[}12{]}“), then $\man{\spa}$ is also
minimum, and is called a \emph{Riemannian Rho Minimum $\sigma$ -representation
(RiRhoMi$\sigma$R) of $\man{\rva}$ .}
\item If $\spa$ is a Minimum $\sigma$R of $\rva$ (Theorem 3 of \cite{Menegaz2015}),
then $\man{\spa}$ is also minimum, and is called a \emph{Riemannian
Minimum $\sigma$-representation (RiMi$\sigma$R) of $\man{\rva}$
.}
\end{itemize}
\end{cor}

The proof of Corollary \ref{cor:Riemannian-sr-minimum-numbers} is
given in Appendix \ref{proof:Corollary_particularRiSRs}.

With Theorem \ref{thm:Euclidean-to-Riemannian-sigma-rep} and Corollary
\ref{cor:Riemannian-sr-minimum-numbers}, we can find an $\risr$
($\weight_{i}^{m}>0$ for every $i=1,\ldots,\nsp$) by first finding
a normalized $\sr$ in the tangent space of the considered manifold;
each normalized $\sr$s (cf. \cite{Menegaz2015} and \cite{Menegaz2016})
have their associated $\risr$s (cf. Corollary \ref{cor:Riemannian-sr-minimum-numbers}).
For instance, suppose we want to calculate the normalized $\mbox{RiMi\ensuremath{\sigma}R}$
of $\man{\rva}\in\man{\rvset}_{\mana}$ (Corollary \ref{cor:Riemannian-sr-minimum-numbers});
that is, we want\footnote{For a set $\man{\xi}:=\{\man{\xi}_{i},\weight_{i}^{m,j},\weight_{i}^{c,j},\weight_{i}^{cc,j}\}$,
if $\weight_{i}^{m,j}=\weight_{i}^{c,j}=\weight_{i}^{cc,j}$ for every
$j=1,...,$l; then we write $\weight_{i}:=\weight_{i}^{m,j}$ and
$\{\man{\xi}_{i},\weight_{i}\}=\man{\xi}$.} 
\[
\mbox{\ensuremath{\man{\spa}}}=\left\{ \man{\spa}_{i},\weight_{i}\right\} _{i=1}^{\dimstate+1}=\mbox{RiMi\ensuremath{\sigma}R}\Big(\est{\man{\state}}_{k-1|k-1},\est{\man{\cov}}_{\man{\state\state}}^{k-1|k-1}\Big).
\]
We can compute the $\mbox{Mi\ensuremath{\sigma}R}$ (Theorem 3 of
\cite{Menegaz2015})
\[
\spa=\big\{\spa_{i},\weight_{i}\big\}_{i=1}^{\dimstate+1}:=\mbox{Mi\ensuremath{\sigma}R}\Big([0]_{\dimstate\times1},\est{\man{\cov}}_{\man{\state\state}}^{k-1|k-1}\Big),
\]
and then, from Theorem \ref{thm:Euclidean-to-Riemannian-sigma-rep},
we would have
\[
\mbox{\ensuremath{\man{\spa}}}=\big\{\rexp_{\est{\man{\state}}_{k-1|k-1}}\spa_{i},\weight_{i}\big\}_{i=1}^{\dimstate+1}.
\]
The work \cite{Hauberg2013} introduced this technique {[}cf. (11)
to (17) therein{]}, and here, with Theorem \ref{thm:Euclidean-to-Riemannian-sigma-rep}
and Corollary \ref{cor:Riemannian-sr-minimum-numbers}, we provide
its formal justification and required assumptions.

\section{Riemannian Unscented Transformations}

\label{sec:Riemannian-Unscented-Transformat}

Essentially, a UT is an approximation of the joint pdf of two functionally-related
random vectors by two weighted sets. For a Riemannian extension of
the UT, we develop likewise.

For the natural numbers $l\geq2$ and $\nsp\geq1$, consider i) a
function $f:\mana\rightarrow\manb$, ii) the random points $\man{\rva}\sim(\mean{\man{\rva}},\man{\moment}_{\man{\rva}}^{2},...,\man{\moment}_{\man{\rva}}^{l}){}_{\mana^{n}}$
and $\man{\rvb}:=\funca(\man{\rva})\sim(\mean{\man{\rvb}},\man{\moment}_{\man{\rvb}}^{2},...,\moment_{\man{\rvb}}^{l})_{\manb^{\eta}}$,
and iii) the sets
\begin{multline*}
\man{\chi}:=\Big\{\man{\spa}_{i},\weight_{i}^{m},\weight_{i}^{c,j},\weight_{i}^{cc,j}:\man{\chi}_{i}\in\mana;\\
j=1,...,l;\,\weight_{i}^{m},\weight_{i}^{c,j},\weight_{i}^{cc,j}\neq0\Big\}_{i=1}^{N}\text{ and}
\end{multline*}
\[
\man{\spb}:=\Big\{\man{\spb}_{i},\weight_{i}^{m},\weight_{i}^{c,j},\weight_{i}^{cc,j}:\man{\spb}_{i}=f(\man{\spa}_{i});\,j=1,...,l\Big\}_{i=1}^{\nsp}.
\]
\begin{defn}[\emph{Ri$l$UT; Definition of 9.2 \cite{Menegaz2016}}]
\label{def:Riemannian-Unscented-Transformation} If $\man{\smean}_{\man{\spa}}=\mean{\man{\rva}}$
and $\man{\smoment}_{\man{\spa}}^{j}=\man{\moment}_{\man{\rva}}^{j}$
for every $j=2,\ldots,l$; then the $l$\emph{th order Riemannian
Unscented Transformation (Ri$l$UT)} is defined by 
\begin{multline*}
\mbox{Ri}l\mbox{UT}:\big(\pfunc,\mean{\man{\rva}},\man{\moment}_{\man{\rva}}^{2},...,\man{\moment}_{\man{\rva}}^{l}\big)\mapsto\\
(\man{\smean}_{\man{\spb}},\man{\smoment}_{\man{\spb}}^{2},...,\man{\smoment}_{\man{\spb}}^{l},\man{\smoment}_{\man{\spa}\man{\spb}}^{2},...,\man{\smoment}_{\man{\spa}\man{\spb}}^{l}).
\end{multline*}
$\man{\spa}$ is called the \emph{independent set} of $\mbox{Ri}l\mbox{UT}$,
and $\man{\spb}$ its \emph{dependent} set.
\end{defn}
Every Ri$l$th$\nsp\sigma$R\emph{ }is an independent set of an $\mbox{Ri}l\mbox{UT}$.
If $l=2$ or $l$ is irrelevant for a given discussion, we can omit
the reference to $l$ and write $\riUT:=\mbox{Ri}2\mbox{UT}$.

$\mbox{Ri}l\mbox{UT}$s are generalizations of $l\UT$s; every $l\UT$
is an $\mbox{Ri}l\mbox{UT}$, and every $\mbox{Ri}l\mbox{UT}$ with
Euclidean points is an $l\UT$. This follows directly from the last
paragraph of Sections \ref{subsec:Statistics-of-random} and of \ref{sec:Statistics-of-weighted}.

An Ri$l$UT can be viewed as a mapping from $2$ random points $\man{\rva}\in\rvset_{\mana}$
and $\man{\rvb}:=\man f(\man X)$ to two Riemannian sets $\man{\spa}$
and $\man{\gamma}$ acting as a \textit{discrete approximation} of
the \textit{\emph{joint}} pdf of $(\man{\rva},\man{\rvb})$. For instance,
an Ri$2$UT can be viewed as the following approximation (this interpretation
is inspired on \cite{Sarkka2007}) 
\[
\left(\begin{array}{c}
\man{\rva}\\
\man{\rvb}
\end{array}\right)\approx\left(\begin{array}{c}
\tilde{\man{\rva}}\\
\tilde{\man{\rvb}}
\end{array}\right)\sim\left(\left(\begin{array}{c}
\man{\smean}_{\man{\spa}}\\
\man{\smean}_{\man{\spb}}
\end{array}\right),\left(\begin{array}{cc}
\man{\scov}_{\man{\spa}\man{\spa}} & \man{\scov}_{\man{\spa}\man{\spb}}\\
\man{\scov}_{\man{\spa}\man{\spb}}^{T} & \man{\scov}_{\man{\spb}\man{\spb}}
\end{array}\right)\right).
\]

\section{Riemannian Unscented Kalman Filters}

\label{sec:Riemannian-Unscented-Filters}

At this point, we still need to develop i) Riemannian systems; and
ii) state correction equations. First, UKFs estimate systems with
random vectors {[}cf. (\ref{eq:additive-system}) and (\ref{eq:general-system}){]};
thus, for Riemannian UKFs (RiUKFs), we define \textit{systems with
Riemannian random points} (Section \ref{subsec:Riemannian-Systems}).
Second, three steps compose UKFs: 1) state prediction, 2) measurement
prediction, and 3) state correction (cf.\cite{Menegaz2015} and \cite{Menegaz2016}).
The Riemannian extensions of steps 1) and 2) are trivial: since UTs
compose steps 1 and 2, we extend them with \textit{$\riUT$s}. But
we still must extend step 3 (Section \ref{subsec:Correction-equations}).
In possession of these two results, we define RiUKFs and provide a
list of some particular forms (Section \ref{subsec:Riemannian-Unscented-Filters}).

\subsection{Riemannian Dynamics Systems}

\label{subsec:Riemannian-Systems}

Up to this point, we have focused on results regarding points on manifolds.
In this section, we focus on results for dynamic state-space systems
on Riemannian manifolds.

The \emph{Riemannian (stochastic discrete-time dynamic) system} in
its \emph{general} form is given by the following pair of equations:
\begin{equation}
\man{\state}_{k}=\pfunc_{k}\left(\man{\state}_{k-1},\man{\pnoise}_{k}\right),\,\man{\meas}_{k}=\mfunc_{k}\left(\man{\state}_{k},\man{\mnoise}_{k}\right)\label{eq:Riemannian-general-system}
\end{equation}
where $k$ is the time step; $\man{\state}_{k}\in\man{\rvset}_{\mana_{\man{\state}}^{\dimstate}}$
the internal state; $\man{\meas}_{k}\in\man{\rvset}_{\mana_{\man{\meas}}^{\dimmeas}}$
is the measured output; $\man{\pnoise}_{k}\in\man{\rvset}_{\mana_{\man{\pnoise}}^{\dimpnoise}}$
the process noise; and $\man{\mnoise}_{k}\in\man{\rvset}_{\mana_{\man{\mnoise}}^{\dimmnoise}}$
the measurement noise. The noises $\man{\pnoise}_{k}$ and $\man{\mnoise}_{k}$
are uncorrelated, $\man{\pnoise}_{k}$ has mean $\mean{\man{\pnoise}}_{k}$
and covariance $\boldsymbol{Q}_{k}$, and $\man{\mnoise}_{k}$ mean
$\mean{\man{\mnoise}}_{k}$ and covariance $\boldsymbol{R}_{k}$.

We also want to consider an\emph{ additive variant of (\ref{eq:Riemannian-general-system})}
because filters for this class of systems are computationally cheaper.
This additive variant of (\ref{eq:Riemannian-general-system}) would
have i) $\man{\pnoise}_{k}$ acting on $\pfunc_{k}(\man{\state}_{k-1})$
by ``adding'' its mean to the mean of $\pfunc_{k}(\man{\state}_{k-1})$
and its covariance to the covariance of $\pfunc_{k}(\man{\state}_{k-1})$,
and ii) $\man{\mnoise}_{k}$ acting similarly on $\mfunc_{k}(\man{\state}_{k})$.
We can work with sums in tangent spaces using the following proposition.

\begin{prop}[Proposition 8.2 of \cite{Menegaz2016}]
\label{prop:Addition of Riemannian random points}Consider a Riemannian
point $\man{\rva}\sim(\bar{\man{\rva}},\man{\cov}_{\man{\rva}\man{\rva}})_{\mana^{n}}$
and a random vector $p\sim(\bar{p},\cov_{pp})_{T_{\bar{\manpta}}\mana^{n}}$.
If $\maxdefdomain(\man{\mean{\rva}})$ is convex, and $\bar{p}\in\gball(\man{\mean{\rva}},r)\cap\cutlocus(\man{\mean{\rva}})$
where $0<r\le\frac{1}{2}\min\{\rinj(\mana),\pi/\sqrt{\kappa})$ and
$\kappa$ is an upper bound of the sectional curvatures of $\mathcal{N}$;
then
\begin{equation}
\rexp_{\mean{\man{\rva}}}\big[\rlogb{\mean{\man{\rva}}\man{\rva}}+p\big]\sim\big(\rexp_{\bar{\man{\rva}}}\bar{p},\man{\cov}_{\man{\rva}\man{\rva}}+\cov_{pp}\big)_{\mana_{\man{\state}}}.\label{eq:Riemannian-addition}
\end{equation}
\end{prop}

The proof of Proposition \ref{prop:Addition of Riemannian random points}
is in Appendix \ref{proof:proposition1}.

Consider this proposition twice: one for the process function with
$\manpta=\pfunc_{k}(\man{\state}_{k-1})$ and $p=\pnoise_{k}$ , and
the other for the measurement function with $\manpta=\mfunc_{k}(\man{\state}_{k})$
and $p=\mnoise_{k}$. Using this reasoning, we define the \emph{additive}
\emph{Riemannian (stochastic discrete-time dynamic) system} as follows
\{equation (9.20) of \cite{Menegaz2016}\}: 
\begin{align}
\man{\state}_{k} & =\mbox{\ensuremath{\rexp}}_{\overline{\pfunc_{k}\big(\man{\state}_{k-1}\big)}}\left[\rlog_{\overline{\pfunc_{k}\big(\man{\state}_{k-1}\big)}}\pfunc_{k}\big(\man{\state}_{k-1}\big)+\pnoise_{k}\right]\nonumber \\
\man{\meas}_{k} & =\mbox{\ensuremath{\rexp}}_{\overline{\mfunc_{k}\big(\man{\state}_{k}\big)}}\left[\rlog_{\overline{\mfunc_{k}\big(\man{\state}_{k}\big)}}\mfunc_{k}\big(\man{\state}_{k}\big)+\mnoise_{k}\right];\label{eq:Riemannian-additive-system}
\end{align}
where $\man{\state}_{k}\in\man{\rvset}_{\mana_{\man{\state}}^{\dimstate}}$
, $\man{\meas}_{k}\in\man{\rvset}_{\mana_{\man{\meas}}^{\dimmeas}}$
, $\pnoise_{k}\in T_{\pfunc_{k}\left(\man{\state}_{k-1}\right)}\mana_{\man{\state}}^{\dimstate}$,
and $\mnoise_{k}\in T_{\pfunc_{k}\left(\man{\state}_{k-1}\right)}\mana_{\man{\meas}}^{\dimmeas}$.
The noise $\pnoise_{k}$ has mean $\mean{\pnoise}_{k}\in T_{\pfunc_{k}\left(\man{\state}_{k-1}\right)}\mana_{\man{\state}}^{\dimstate}$
and covariance $\covpnoise_{k}\in T_{\pfunc_{k}\left(\man{\state}_{k-1}\right)}\mana_{\man{\state}}^{\dimstate}\times T_{\pfunc_{k}\left(\man{\state}_{k-1}\right)}\mana_{\man{\state}}^{\dimstate}$,
and $\mnoise_{k}$ mean $\mean{\mnoise}_{k}\in T_{\pfunc_{k}\left(\man{\state}_{k-1}\right)}\mana_{\man{\meas}}^{\dimmeas}$
and covariance $\covmnoise_{k}\in T_{\pfunc_{k}\left(\man{\state}_{k-1}\right)}\mana_{\man{\meas}}^{\dimmeas}\times T_{\pfunc_{k}\left(\man{\state}_{k-1}\right)}\mana_{\man{\meas}}^{\dimmeas}$.
Note that $\pnoise_{k}$ is defined in the tangent space $T_{\pfunc_{k}\left(\man{\state}_{k-1}\right)}\mana_{\man{\state}}^{\dimstate}$
and $\mnoise_{k}$ in $T_{\pfunc_{k}\left(\man{\state}_{k-1}\right)}\mana_{\man{\meas}}^{\dimmeas}$.
In Remark \ref{rem:Alternative-additive-Riemannian-system}, we discuss
an alternative definition in which these noises belong to Riemannian
manifolds. An example with the unit sphere manifold of dimension 3,
$\sset 3$, is provided in Section \ref{sec:Numerical-Example}.

To the best of our knowledge, (\ref{eq:Riemannian-additive-system})
is the \emph{first consistent additive-noise Riemannian system}. Although
the literature has introduced additive-noise discrete-time UKFs for
some Riemannian manifolds, we could not find any additive-noise system
retaining the random point in the working manifolds; even for simple
manifolds such as $\sset 3$ (cf. \cite{Crassidis2003,Chang2016,Vartiainen2014}).

If $\mana_{\man{\state}}^{\dimstate}=\realset^{\dimstate}$ and $\mana_{\man{\meas}}^{\dimmeas}=\realset^{\dimmeas}$
then (\ref{eq:Riemannian-additive-system}) is the \emph{additive
system} (\ref{eq:additive-system}). This is a direct consequence
of the following results: for $a,b\in\realset^{n}$ $\rlog_{a}b=b-a$
and $\mbox{\ensuremath{\rexp}}_{a}b=b+a$. 

Sometimes, only one of the two equations in (\ref{eq:Riemannian-general-system})
can be written with additive-noise as in (\ref{eq:Riemannian-additive-system}).
In this case, we define the following two partially-additive Riemannian
systems:
\begin{align}
\man{\state}_{k} & =\pfunc_{k}\left(\man{\state}_{k-1},\man{\pnoise}_{k}\right)\nonumber \\
\man{\meas}_{k} & =\mbox{\ensuremath{\rexp}}_{\overline{\mfunc_{k}\big(\man{\state}_{k}\big)}}\left[\rlog_{\overline{\mfunc_{k}\big(\man{\state}_{k}\big)}}\mfunc_{k}\big(\man{\state}_{k}\big)+\mnoise_{k}\right];\label{eq:Riemannian-partiallly-additive-system1}
\end{align}
and
\begin{align}
\man{\state}_{k} & =\mbox{\ensuremath{\rexp}}_{\overline{\pfunc_{k}\big(\man{\state}_{k-1}\big)}}\left[\rlog_{\overline{\pfunc_{k}\big(\man{\state}_{k-1}\big)}}\pfunc_{k}\big(\man{\state}_{k-1}\big)+\pnoise_{k}\right]\nonumber \\
\man{\meas}_{k} & =\mfunc_{k}\left(\man{\state}_{k},\man{\mnoise}_{k}\right).\label{eq:Riemannian-partiallly-additive-system2}
\end{align}

\begin{rem}[]
\label{rem:Alternative-additive-Riemannian-system}System (\ref{eq:Riemannian-additive-system})
is defined with tangent space process and measurement noises. An alternative
definition in which these noises belong to Riemannian manifolds is
the following:
\begin{align*}
\man{\state}_{k} & =\mbox{\ensuremath{\rexp}}_{\overline{\pfunc_{k}\left(\man{\state}_{k-1}\right)}}\left[\rlog_{\overline{\pfunc_{k}\left(\man{\state}_{k-1}\right)}}\pfunc_{k}\left(\man{\state}_{k-1}\right)+\rlog_{\overline{\pfunc_{k}\left(\man{\state}_{k-1}\right)}}\man{\pnoise}_{k}\right]\\
\man{\meas}_{k} & =\mbox{\ensuremath{\rexp}}_{\overline{\mfunc_{k}\left(\man{\state}_{k}\right)}}\left[\rlog_{\overline{\mfunc_{k}\left(\man{\state}_{k}\right)}}\mfunc_{k}\left(\man{\state}_{k}\right)+\rlog_{\overline{\mfunc_{k}\left(\man{\state}_{k}\right)}}\man{\mnoise}_{k}\right];
\end{align*}
where $\man{\state}_{k}\in\man{\rvset}_{\mana_{\man{\state}}^{\dimstate}}$
, $\man{\meas}_{k}\in\man{\rvset}_{\mana_{\man{\meas}}^{\dimmeas}}$
, $\man{\pnoise}_{k}\in\man{\rvset}_{\mana_{\man{\state}}^{\dimstate}}$,
and $\man{\mnoise}_{k}\in\man{\rvset}_{\mana_{\man{\meas}}^{\dimmeas}}$.
In this case, it would be interesting to assume one of the following
two cases:

\begin{enumerate}
\item That are known i) the means of $\man{\pnoise}_{k}$ and $\man{\mnoise}_{k}$—e.g.,
$\mean{\man{\pnoise}}_{k}\in\mana_{\man{\state}}^{\dimstate}-\cutlocus(\overline{\pfunc_{k}\left(\man{\state}_{k-1}\right)})$
and $\mean{\man{\mnoise}}_{k}\in\mana_{\man{\meas}}^{\dimmeas}-\cutlocus(\overline{\mfunc_{k}\left(\man{\state}_{k}\right)})$—,
b) the covariance of $\man{\pnoise}_{k}$ with respect to $\mean{\man{\pnoise}}_{k}$
at $\overline{\pfunc_{k}\left(\man{\state}_{k-1}\right)}$, and iii)
the covariance of $\man{\mnoise}_{k}$ with respect to $\mean{\man{\mnoise}}_{k}$
at $\overline{\mfunc_{k}\left(\man{\state}_{k}\right)}$.
\item That the means and covariances of $\rlog_{\overline{\pfunc_{k}\left(\man{\state}_{k-1}\right)}}\man{\pnoise}_{k}$
and $\rlog_{\overline{\mfunc_{k}\left(\man{\state}_{k}\right)}}\man{\mnoise}_{k}$
are known—e.g., the means $\mean{\pnoise}_{k}\in T_{\overline{\pfunc_{k}\left(\man{\state}_{k-1}\right)}}\man{\rvset}_{\mana_{\man{\state}}^{\dimstate}}$
and $\mean r_{k}\in T_{\overline{\mfunc_{k}\left(\man{\state}_{k}\right)}}\man{\rvset}_{\mana_{\man{\meas}}^{\dimmeas}}$;
and the covariances $\covpnoise_{k}\in T_{\overline{\pfunc_{k}\left(\man{\state}_{k-1}\right)}}\man{\rvset}_{\mana_{\man{\state}}^{\dimstate}}\times T_{\overline{\pfunc_{k}\left(\man{\state}_{k-1}\right)}}\man{\rvset}_{\mana_{\man{\state}}^{\dimstate}}$
and $\covmnoise_{k}\in T_{\overline{\mfunc_{k}\left(\man{\state}_{k}\right)}}\man{\rvset}_{\mana_{\man{\state}}^{\dimstate}}\times T_{\overline{\mfunc_{k}\left(\man{\state}_{k}\right)}}\man{\rvset}_{\mana_{\man{\meas}}^{\dimmeas}}$.
\end{enumerate}
\end{rem}

\subsection{Correction equations}

\label{subsec:Correction-equations}

In this section, we introduce Riemannian extensions of the UKFs correction
equations. Finding these extensions is not trivial because their Euclidean
versions include vector operations (cf. \cite{Menegaz2015}), which
are not defined for all Riemannian manifolds. Thus, we proceed by
first considering the simpler case $\man{\rvset}_{\mana_{\man{\state}}^{\dimstate}}=\man{\rvset}_{\mana_{\man{\meas}}^{\dimmeas}}$.

\subsubsection{State and measurement in the same manifold}

\label{subsec:State-and-measurement}

Suppose that $\mana_{\man{\state}}^{\dimstate}=\mana_{\man{\meas}}^{\dimmeas}$
and the measurements $\outcome{\man{\meas}}_{1}$, ..., $\outcome{\man{\meas}}_{k}$
have been acquired. Define the following random points\footnote{For the random points $\man{\rva}$ and $\man{\rvb}$ and the outcomes
$\outcome{\man{\rvb}}_{1}$, ..., $\outcome{\man{\rvb}}_{l}$ of $\man{\rvb}$;
the random point $\man{\rva}$|$\outcome{\man{\rvb}}_{1:k-1}$ stands
for $\man{\rva}$ conditioned to $\man{\rvb}_{i}=\outcome{\man{\rvb}}_{i}$
for every $i=1,$..., $l$.}
\begin{align*}
\man{\state}_{k|k-1} & :=\man{\state}_{k}|\outcome{\man{\meas}}_{1:k-1}\\
\man{\state}_{k|k} & :=\man{\state}_{k}|\outcome{\man{\meas}}_{1:k}\\
\man{\meas}_{k|k-1} & :=\man{\meas}_{k}|\outcome{\man{\meas}}_{1:k-1},
\end{align*}
and the following projections on the tangent space of $\man{\state}_{k|k-1}$
\begin{align}
\state_{k|k-1}^{TM} & :=\rlog_{\mean{\man{\state}}_{k|k-1}}\man{\state}_{k|k-1}\label{eq:corrected-state-in-tangent-space}\\
\state_{k|k}^{TM} & :=\rlog_{\mean{\man{\state}}_{k|k-1}}\man{\state}_{k|k}\nonumber \\
\meas_{k|k-1}^{TM} & :=\rlog_{\mean{\man{\state}}_{k|k-1}}\man{\meas}_{k|k-1}\label{eq:predicted-measurement-in-tangent-space}\\
\outcome{\meas}_{k}^{TM} & :=\rlog_{\mean{\man{\state}}_{k|k-1}}\outcome{\man{\meas}}_{k}.\label{eq:realization-of-measurement-in-tangent-space}
\end{align}
Let i) $\man{\state}_{k|k-1}$ and $\man{\meas}_{k|k-1}$ be characterized
by their projection on the tangent space of $\man{\state}_{k|k-1}$
according to the following equation: 
\begin{multline}
\big[\begin{array}{cc}
\state_{k|k-1}^{TM} & \meas_{k|k-1}^{TM}\big]^{T}\end{array}\sim\\
\normal\left(\left[\begin{array}{c}
[0]_{\dimstate,1}\\
\meas_{k|k-1}^{TM}
\end{array}\right],\left[\begin{array}{cc}
\man{\cov}_{\man{\state\state}}^{k|k-1} & \man{\cov}_{\man{\state\meas}}^{k|k-1}\\
\left(\man{\cov}_{\man{\state y}}^{k|k-1}\right)^{T} & \man{\cov}_{\man{\meas\meas}}^{k|k-1}
\end{array}\right]\right);\label{eq:predicted-tangential-statistics-are-normal}
\end{multline}
and ii) the projection $\state_{k|k}^{TM}$ be given by the following
linear correction of $\state_{k|k-1}^{TM}$ 
\begin{equation}
\state_{k|k}^{TM}=\state_{k|k-1}^{TM}+\man{\kgain}_{k}\left(\outcome{\meas}_{k}^{TM}-\meas_{k|k-1}^{TM}\right),\label{eq:correction-equation-of-xk-tangent}
\end{equation}
where $\man{\kgain}_{k}\in\realset^{\dimstate\times\dimstate}$ is
a gain matrix. From known results of the Kalman filtering theory (cf.
\cite{Anderson1979}), we have
\begin{equation}
\man{\kgain}_{k}:=\man{\cov}_{\man{\state\meas}}^{k|k-1}\left(\man{\cov}_{\man{\meas\meas}}^{k|k-1}\right)^{-1},\label{eq:Riemannian-kalman-gain}
\end{equation}
and $\state_{k|k}^{TM}\sim\normal(\mean{\state}_{k|k}^{TM},\man{\cov}_{\man{\state\state}}^{k|k-1,\mean{\man{\state}}_{k|k-1}})$
where
\begin{align}
\mean{\state}_{k|k}^{TM} & :=\man{\kgain}_{k}\left(\outcome{\meas}_{k}^{TM}-\mean{\meas}_{k|k-1}^{TM}\right)\label{eq:Riemannian-corrected-mean}\\
\man{\cov}_{\man{\state\state}}^{k|k,\mean{\man{\state}}_{k|k-1}} & :=\man{\cov}_{\man{\state\state}}^{k|k-1}-\left(\man{\kgain}_{k}\right)\man{\cov}_{\man{\meas\meas}}^{k|k-1}\left(\man{\kgain}_{k}\right)^{T}.\label{eq:Riemannian-corrected-cov-in-tangent-of-xk_k-1}
\end{align}
From (\ref{eq:corrected-state-in-tangent-space}), we have
\begin{equation}
\mean{\man{\state}}_{k|k}=\rexp_{\mean{\man{\state}}_{k|k-1}}\state_{k|k}^{TM}.\label{eq:Riemannian-corrected-mean-1}
\end{equation}

The matrix $\man{\cov}_{\man{\state\state}}^{k|k-1,\mean{\man{\state}}_{k|k-1}}$
is the covariance of $\man{\state}_{k|k}$ relative to $\mean{\man{\state}}_{k|k}$
\emph{at} $\mean{\man{\state}}_{k|k-1}$. We want the covariance $\man{\cov}_{\man{\state\state}}^{k|k}:=\man{\cov}_{\man{\state\state}}^{k|k,\mean{\man{\state}}_{k|k}}$
of $\man{\state}_{k|k}$ \emph{at} $\mean{\man{\state}}_{k|k}$, and
the following theorem from \cite{Hauberg2013} provides the mechanism
to obtain $\man{\cov}_{\man{\state\state}}^{k|k,\mean{\man{\state}}_{k|k}}$
from $\man{\cov}_{\man{\state\state}}^{k|k,\mean{\man{\state}}_{k|k-1}}$. 

\begin{thm}[Parallel Transport of a Bilinear Mapping \cite{Hauberg2013}]
\label{thm:Hauberg-matrix-parellel-transport}Let $\cov$ be a \emph{symmetric}
\emph{bilinear} mapping on the tangent space $T_{\manpta}\mana$ of
the Riemannian manifold $\mana$ at $\manpta\in\mana$, and $\curve:[0,1]\rightarrow\mana$
a differentiable curve on $\mana$ with $\curve(0)=\manpta$. Since
$\cov$ is symmetric, it can be written as
\[
\cov=\sum_{i=1}^{n}\lambda_{i}v_{i}v_{i}^{T}
\]
where ($v_{1},$ ..., $v_{n}$) is an orthonormal basis of $T_{\manpta}\mana$,
and each $\lambda_{i}$ is the eigenvalue of $\cov$ associated with
the eigenvector $v_{i}$. Let $v_{i}(t)$ be the parallel transport
of $v_{i}$ along $\curve(t)$ (Definition \ref{def:parallel-transport}).
With this,
\begin{equation}
\cov_{t}:=\sum_{i=1}^{n}\lambda_{i}v_{i}(t)v_{i}(t)^{T}\label{eq:parallel-transport-theorem}
\end{equation}
is the \emph{parallel transport of $\cov$} along $\curve(t)$. 
\end{thm}

When we do not know the closed form of a tangent vector parallel transport,
we can use a numerical approach such as the \emph{Schild's Ladder}
(cf. \cite{Hauberg2013}; see \cite{Lorenzi2013} for other implementations
and algorithms of parallel transports).

We obtain $\man{\cov}_{\man{\state\state}}^{k|k}$ by performing the
parallel transport of $\man{\cov}_{\man{\state\state}}^{k|k,\mean{\man{\state}}_{k|k}}$from
$\mean{\man{\state}}_{k|k-1}$ to $\mean{\man{\state}}_{k|k}$ as
follows:
\begin{equation}
\man{\cov}_{\man{\state\state}}^{k|k}=\partransp\Big(\man{\cov}_{\man{\state\state}}^{k|k,\mean{\man{\state}}_{k|k}},\mean{\man{\state}}_{k|k-1},\mean{\man{\state}}_{k|k}\Big),\label{eq:Riemannian-corrected-cov}
\end{equation}
where 
\begin{eqnarray*}
\partransp: & \textrm{Sym}\left(T_{\manpta}\mana\right)\times\mana\times\mana & \rightarrow\textrm{Sym}\left(T_{\manptb}\mana\right)\\
 & \left(\cov^{\manpta},\manpta,\manptb\right) & \mapsto\cov^{\manptb}
\end{eqnarray*}
is the function mapping $\textrm{Sym}(T_{\manpta}\mana)\times\mana\times\mana$
to $\textrm{Sym}(T_{\manptb}\mana)$ according to (\ref{eq:parallel-transport-theorem}),
and $\textrm{Sym}(T_{\manpta}\mana)$ denotes the space of symmetric
matrices of $T_{\manpta}\mana$.

With this, we can define a UKF for Riemannian systems when $\mana_{\man{\state}}^{\dimstate}=\mana_{\man{\meas}}^{\dimmeas}$.
Let us now consider the original more general case.

\subsubsection{State and measurement in different manifolds}

\label{subsec:State-and-measurement=00005B}

If $\man{\state}_{k}$ belongs to a manifold $\man{\rvset}_{\mana_{\man{\state}}^{\dimstate}}$
and $\man{\meas}_{k}$ to another manifold $\man{\rvset}_{\mana_{\man{\meas}}^{\dimmeas}}$,
then we can not define $\meas_{k|k-1}^{TM}$ as in (\ref{eq:predicted-measurement-in-tangent-space})
and $\outcome{\meas}_{k}^{TM}$ as in (\ref{eq:realization-of-measurement-in-tangent-space});
consequently, neither $\state_{k|k}^{TM}$ as in (\ref{eq:correction-equation-of-xk-tangent}).

Since we know the correction equations when $\mana_{\man{\state}}^{\dimstate}=\mana_{\man{\meas}}^{\dimmeas}$,
we can look for a manifold of which both $\mana_{\man{\state}}$ and
$\mana_{\man{\meas}}$ are submanifolds. The simplest of such a class
is $\mana_{\man{\state}}\times\mana_{\man{\meas}}$—the Cartesian
product of two Riemannian manifolds is a Riemannian manifold \cite{DoCarmo1992}.

Suppose $\state_{k|k-1}^{TM}$ and $\meas_{k|k-1}^{TM}$ are jointly
Gaussian random vectors according to (\ref{eq:predicted-tangential-statistics-are-normal}).
Define i) the Riemannian Manifold $\mana_{\man{\state},\man{\meas}}:=\mana_{\man{\state}}\times\mana_{\man{\meas}}$;
ii) the points $\man c:=(\man c_{\man{\state}},\man c_{\man{\meas}})\in\mana_{\man{\state},\man{\meas}}$,
$\man b_{\man{\state}}\in\mana_{\man{\state}}$, and $\man b_{\man{\meas}}\in\mana_{\man{\meas}}$
(these points are chosen); and the following random vector belonging
to $T_{\man c}\mana_{\man{\state},\man{\meas}}$: 
\begin{multline*}
\state_{k|k,**}^{T_{\man c}\mana_{\man{\state},\man{\meas}}}:=\rlog_{\man c}\big[\man{\state}_{k|k-1},\,\man b_{\man{\meas}}\big]^{T}\\
+\man{\kgain}_{k,**}\big(\rlog_{\man c}[\man b_{\man{\state}},\,\outcome{\man{\meas}}_{k}]^{T}-\rlog_{\man c}[\man b_{\man{\state}},\man{\meas}_{k|k-1}]^{T}\big)
\end{multline*}
where $\man{\kgain}_{k,**}\in\realset^{\left(\dimstate+\dimmeas\right)\times\left(\dimstate+\dimmeas\right)}$
is a gain matrix. The tangent vector $\state_{k|k,**}^{T_{\man c}\mana_{\man{\state},\man{\meas}}}$
is clearly related with $\state_{k|k}^{TM}$ by
\begin{equation}
\state_{k|k}^{TM}:=\left[\hat{\state}_{k|k,**}^{T_{\man c}\mana_{\man{\state},\man{\meas}}}\right]_{1:\dimstate,1}.\label{eq:augmented-corrected-state-to-non-augmented}
\end{equation}
By finding the mean and covariance of $\state_{k|k,**}^{T_{\man c}\mana_{\man{\state},\man{\meas}}}$,
we find the mean and covariance of $\state_{k|k}^{TM}$.

Since $\state_{k|k-1}^{TM}$ and $\meas_{k|k-1}^{TM}$ are jointly
Gaussian random vectors, it follows that—we use the same reasoning
used to obtain (\ref{eq:Riemannian-kalman-gain}), (\ref{eq:Riemannian-corrected-mean}),
(\ref{eq:Riemannian-corrected-cov-in-tangent-of-xk_k-1}), (\ref{eq:Riemannian-corrected-mean-1}),
and (\ref{eq:Riemannian-corrected-cov})—
\begin{align*}
\man{\cov}_{\man{\state}\man{\state},**}^{k|k-1} & :=\man{\ev}\left\{ \left(\rlog_{\man c}\left[\begin{array}{c}
\man x\\
\man b_{\man{\meas}}
\end{array}\right]-\rlog_{\man c}\left[\begin{array}{c}
\mean{\man{\state}}_{k|k-1}\\
\man b_{\man{\meas}}
\end{array}\right]\right)\left(\diamond\right)^{T}\right\} \\
\man{\cov}_{\man{\meas}\man{\meas},**}^{k|k-1} & :=\man{\ev}\left\{ \left(\rlog_{\man c}\left[\begin{array}{c}
\man b_{\man{\state}}\\
\man{\meas}
\end{array}\right]-\rlog_{\man c}\left[\begin{array}{c}
\man b_{\man{\state}}\\
\mean{\man{\meas}}_{k|k-1}
\end{array}\right]\right)\left(\diamond\right)^{T}\right\} 
\end{align*}
\begin{multline*}
\man{\cov}_{\man{\state}\man{\meas},**}^{k|k-1}:=\man{\ev}\left\{ \left(\rlog_{\man c}\left[\begin{array}{c}
\man x\\
\man b_{\man{\meas}}
\end{array}\right]-\rlog_{\man c}\left[\begin{array}{c}
\mean{\man{\state}}_{k|k-1}\\
\man b_{\man{\meas}}
\end{array}\right]\right)\right.\\
\times\left.\left(\rlog_{\man c}\left[\begin{array}{c}
\man b_{\man{\state}}\\
\man{\meas}
\end{array}\right]-\rlog_{\man c}\left[\begin{array}{c}
\man b_{\man{\state}}\\
\mean{\man{\meas}}_{k|k-1}
\end{array}\right]\right)^{T}\right\} ;
\end{multline*}
thus, the mean and covariance of $\state_{k|k,**}^{T_{\man c}\mana_{\man{\state},\man{\meas}}}$
are given by
\begin{align}
\man{\kgain}_{k,**}:= & \man{\cov}_{\man{\state}\man{\meas},**}^{k|k-1}\diag\Big([0]_{\dimstate\times\dimstate},\nonumber \\
 & \man{\ev}\{(\rlog_{\man c_{\man{\meas}}}\man{\meas}-\rlog_{\man c_{\man{\meas}}}\mean{\man{\meas}}_{k|k-1})\left(\diamond\right)^{T}\}^{-1}\Big)\label{eq:UKF-Riemannian-CartesianSolution-3}\\
\mean{\state}_{k|k,**}^{T_{\man c}\mana_{\man{\state},\man{\meas}}}:= & \rlog_{\man c}\left[\begin{array}{c}
\mean{\man{\state}}_{k|k-1}\\
\man b_{\man{\meas}}
\end{array}\right]+\man{\kgain}_{k,**}\rlog_{\man c}\left[\begin{array}{c}
\man b_{\man{\state}}\\
\mean{\man{\meas}}_{k|k-1}
\end{array}\right]\label{eq:UKF-Riemannian-CartesianSolution-4}\\
\man{\cov}_{\man{\state}\man{\state},**}^{k|k,T_{\man c}M}:= & \man{\cov}_{\man{\state}\man{\state},**}^{k|k-1}-\left(\man{\kgain}_{k,**}\right)\man{\cov}_{\man{\meas}\man{\meas},**}^{k|k-1}\left(\man{\kgain}_{k,**}\right)^{T}.\label{eq:UKF-Riemannian-CartesianSolution-6}
\end{align}
We can choose $\man c$, $\man b_{\man{\state}}$ and $\man b_{\man{\meas}}$
arbitrarily, and a particular choice yields the desired correction
equations.

\begin{thm}[Theorem 9.3 of \cite{Menegaz2016}]
\label{thm:General-to-Hauberg-UKF}Given (\ref{eq:augmented-corrected-state-to-non-augmented}),
(\ref{eq:UKF-Riemannian-CartesianSolution-3}), (\ref{eq:UKF-Riemannian-CartesianSolution-4}),
and (\ref{eq:UKF-Riemannian-CartesianSolution-6}); if $\man c_{\man{\state}}=\man b_{\man{\state}}=\est{\man{\state}}_{k|k-1}$
and $\man c_{\man{\meas}}=\man b_{\man{\meas}}=\est{\man{\meas}}_{k|k-1}$,
then 
\begin{equation}
\state_{k|k}^{TM}=\man{\kgain}_{k}\rlog_{\est{\man{\meas}}_{k|k-1}}\left(\man{\meas}_{k}\right)\label{eq:corrected-state-tangent-final-equation}
\end{equation}
and 
\begin{equation}
\man{\cov}_{\man{\state}\man{\state}}^{k|k,\mean{\man{\state}}_{k|k-1}}=\est{\man{\cov}}_{\man{\state}\man{\state}}^{k|k-1}-\man{\kgain}_{k}\left(\est{\man{\cov}}_{\man{\meas}\man{\meas}}^{k|k-1}\right)^{-1}\left(\man{\kgain}_{k}\right)^{T},\label{eq:corrected-coov-tangent-xk-k-1-final-equation}
\end{equation}
where
\[
\man{\kgain}_{k}:=\man{\cov}_{\man{\state\meas}}^{k|k-1}\left(\man{\cov}_{\man{\meas\meas}}^{k|k-1}\right)^{-1}.
\]
\end{thm}

The proof of Theorem \ref{thm:General-to-Hauberg-UKF} is in Appendix
\ref{proof:TheoremRiUKFcorrection}.

According to this theorem, the correction equations—(\ref{eq:Riemannian-kalman-gain}),
(\ref{eq:Riemannian-corrected-mean}), (\ref{eq:Riemannian-corrected-cov-in-tangent-of-xk_k-1}),
(\ref{eq:Riemannian-corrected-mean-1}), and (\ref{eq:Riemannian-corrected-cov})—are
true even when the state and the measurement belong to different manifolds.
Therefore, we do not have to perform calculations on the bigger manifold
$\mana_{\man{\state},\man{\meas}}$ to calculate $\state_{k|k}^{TM}$
and $\man{\cov}_{\man{\state}\man{\state}}^{k|k,\mean{\man{\state}}_{k|k-1}}$.
Instead, they can be calculated by (\ref{eq:corrected-state-tangent-final-equation})
and (\ref{eq:corrected-coov-tangent-xk-k-1-final-equation}) even
when $\mana_{\man{\state}}\neq\mana_{\man{\meas}}$.

\subsection{New Riemannian Unscented Kalman Filters}

\label{subsec:Riemannian-Unscented-Filters}

At this point, we are endowed with the necessary results to provide
Riemannian extensions of UKFs. At every step time, the final estimates
$\est{\man{\state}}_{k|k}$ and $\est{\man{\cov}}_{\man{\state\state}}^{k|k}$
can be calculated by (\ref{eq:Riemannian-corrected-cov}) and Theorem
\ref{thm:General-to-Hauberg-UKF}. From (\ref{eq:Riemannian-kalman-gain}),
(\ref{eq:Riemannian-corrected-cov}), (\ref{eq:corrected-state-tangent-final-equation}),
and (\ref{eq:corrected-coov-tangent-xk-k-1-final-equation}) these
final estimates require $\est{\man{\state}}_{k|k-1}$, $\est{\man{\cov}}_{\man{\state\state}}^{k|k-1}$,
$\est{\man{\meas}}_{k|k-1}$, $\est{\man{\cov}}_{\man{yy}}^{k|k-1}$,
and $\est{\man{\cov}}_{\man{\state y}}^{k|k-1}$. These last estimates
can be calculated by realizing $\riUT$s in systems (\ref{eq:Riemannian-general-system})
and (\ref{eq:Riemannian-additive-system}). For instance, from (\ref{eq:Riemannian-additive-system}),
Definition \ref{def:Riemannian-Unscented-Transformation} and Proposition
\ref{prop:Addition of Riemannian random points}, the estimates $\est{\man{\state}}_{k|k-1}$,
$\est{\man{\cov}}_{\man{\state\state}}^{k|k-1}$ can be calculated
by
\begin{multline*}
\Big(\est{\man{\meas}}_{k|k-1}^{*},\est{\man{\cov}}_{\man{yy},*}^{k|k-1},\est{\man{\cov}}_{\man{\state y}}^{k|k-1}\Big):=\\
\riUT_{2}\left(h_{k},\est{\man{\state}}_{k|k-1},\est{\man{\cov}}_{\man{\state\state}}^{k|k-1}\right)
\end{multline*}
\begin{align*}
\est{\man{\meas}}_{k|k-1} & :=\rexp_{\est{\man{\meas}}_{k|k-1}^{*}}\mean{\mnoise}_{k}\\
\est{\man{\cov}}_{\man{yy}}^{k|k-1} & :=\est{\man{\cov}}_{\man{yy},*}^{k|k-1}+\covmnoise_{k}.
\end{align*}
By similar formulas, we can obtain $\est{\man{\state}}_{k|k-1}$,
$\est{\man{\cov}}_{\man{\state\state}}^{k|k-1}$, $\est{\man{\meas}}_{k|k-1}$,
$\est{\man{\cov}}_{\man{yy}}^{k|k-1}$, and $\est{\man{\cov}}_{\man{\state y}}^{k|k-1}$
for both (\ref{eq:Riemannian-general-system}) and (\ref{eq:Riemannian-additive-system}).

Below, we introduce the Riemannian UKFs (RiUKFs): UKFs for the Riemannian
systems (\ref{eq:Riemannian-general-system}) and (\ref{eq:Riemannian-additive-system}).
For the filter of (\ref{eq:Riemannian-general-system}), define the
augmented functions
\begin{align}
\pfunc_{k}^{aug}\big([\man{\state}_{k-1},\,\man{\pnoise}_{k}]^{T}\big) & :=\pfunc_{k}\big(\man x_{k-1},\man{\pnoise}_{k}\big)\label{eq:Riemannian-augmented-functions}\\
\mfunc_{k}^{aug}\big([\man{\state}_{k},\,\man{\mnoise}_{k}]^{T}\big) & :=\mfunc_{k}\big(\man x_{k},\man{\mnoise}_{k}\big).\nonumber 
\end{align}

Consider system (\ref{eq:Riemannian-general-system}) and suppose
that i) the initial state is $\man{\state}_{0}\sim\left(\mean{\man{\state}}_{0},\man{\cov}_{\man{\state}\man{\state}}^{0}\right)_{\mana_{\man{\state}}},$
and ii) the measurements $\outcome{\man{\meas}}_{1}$, $\outcome{\man{\meas}}_{2}$,
..., $\outcome{\man{\meas}}_{k_{f}}$ are given. Then the \emph{Riemannian
Augmented Unscented Kalman Filter }(RiAuUKF)\emph{ }is given by the
following algorithm:
\begin{lyxalgorithm}[RiAuUKF; Algorithm 19 of \cite{Menegaz2016}]
\label{alg:RiAuUKF}Set the initial estimates $\est{\man{\state}}_{0|0}:=\mean{\man{\state}}_{0}$
and $\est{\man{\cov}}_{\state\state}^{0|0}:=\man{\cov}_{\state\state}^{0}$.
For $k=1,...,k_{f}$, perform the following steps:
\begin{enumerate}[labelsep=0.1cm,leftmargin=0.45cm]
\item \label{enu:RiAuUKF-State-prediction}State prediction.\textup{
\begin{align*}
\est{\man{\state}}_{k-1|k-1}^{aug} & :=\left[\est{\man{\state}}_{k-1|k-1}^{T},\mean{\man{\pnoise}}_{k}^{T}\right]^{T}\\
\est{\man{\cov}}_{\man{\state\state},aug}^{k-1|k-1} & :=\mbox{\ensuremath{\diag}}\left(\est{\man{\cov}}_{\man{\state\state}}^{k-1|k-1},\man{\covpnoise}_{k}\right)
\end{align*}
\begin{multline}
\Big(\est{\man{\state}}_{k|k-1},\est{\man{\cov}}_{\man{\state\state}}^{k|k-1}\Big):=\\
\riUT_{1}\left(f_{k}^{aug},\est{\man{\state}}_{k-1|k-1}^{aug},\est{\man{\cov}}_{\man{\state\state},aug}^{k-1|k-1}\right).\label{RiAuUKF-RiUT1}
\end{multline}
}
\item \label{enu:RiAuUKF-Measurement-prediction}Measurement prediction.\textup{
\begin{align*}
\est{\man{\state}}_{k|k-1}^{aug} & :=\left[\est{\man{\state}}_{k|k-1}^{T},\mean{\man{\mnoise}}_{k}^{T}\right]^{T}\\
\est{\man{\cov}}_{\man{\state\state},aug}^{k|k-1} & :=\mathcal{\diag}\left(\est{\man{\cov}}_{\man{\state\state}}^{k|k-1},\man{\covmnoise}_{k}\right).
\end{align*}
\begin{multline}
\Big(\est{\man{\meas}}_{k|k-1},\est{\man{\cov}}_{\man{yy}}^{k|k-1},\est{\man{\cov}}_{\man{\state y},a}^{k|k-1}\Big):=\\
\riUT_{2}\left(h_{k}^{aug},\est{\man{\state}}_{k|k-1}^{aug},\est{\man{\cov}}_{\man{\state\state},aug}^{k|k-1}\right)\label{RiAuUKF-RiUT2}
\end{multline}
\[
\est{\man{\cov}}_{\man{\state y}}^{k|k-1}:=\left[\est{\man{\cov}}_{\man{\state y},aug}^{k|k-1}\right]_{\left(1:\dimstate\right),\left(1:\dimmeas\right)}.
\]
}
\item \label{enu:RiAuUKF-State-correction}State correction.\textup{
\begin{align}
\man{\kgain}_{k} & :=\est{\man{\cov}}_{\man{\state}\man{\meas}}^{k|k-1}\Big(\est{\man{\cov}}_{\man{\meas}\man{\meas}}^{k|k-1}\Big)^{-1}\label{eq:RiAuUKF-KalmanGain}\\
\hat{\state}_{k|k}^{TM} & :=\man{\kgain}_{k}\rlog_{\est{\man{\meas}}_{k|k-1}}\big(\outcome{\man{\meas}}_{k}\big)\nonumber \\
\est{\man{\state}}_{k|k} & :=\rexp_{\est{\man{\state}}_{k|k-1}}\big(\hat{\state}_{k|k}^{TM}\big)\nonumber \\
\est{\man{\cov}}_{\man{\state}\man{\state}}^{k|k,\est{\man{\state}}_{k|k-1}} & :=\est{\man{\cov}}_{\man{\state}\man{\state}}^{k|k-1}-\man{\kgain}_{k}\est{\man{\cov}}_{\man{\meas}\man{\meas}}^{k|k-1}\man{\kgain}_{k}^{T}\nonumber \\
\est{\man{\cov}}_{\man{\state}\man{\state}}^{k|k} & :=\partransp\Big(\est{\man{\cov}}_{\man{\state}\man{\state}}^{k|k,\est{\man{\state}}_{k|k-1}},\est{\man{\state}}_{k|k-1},\est{\man{\state}}_{k|k}\Big).\nonumber 
\end{align}
}
\end{enumerate}
\end{lyxalgorithm}
Consider the system (\ref{eq:Riemannian-additive-system}) and suppose
that i) the initial state is $\man{\state}_{0}\sim\left(\mean{\man{\state}}_{0},\man{\cov}_{\man{\state}\man{\state}}^{0}\right)_{\mana_{\man{\state}}},$
and ii) the measurements $\outcome{\man{\meas}}_{1}$, $\outcome{\man{\meas}}_{2}$,
..., $\outcome{\man{\meas}}_{k_{f}}$ are given. Then the \emph{Riemannian
Additive Unscented Kalman Filter }(RiAdUKF) is given by the following
algorithm:
\begin{lyxalgorithm}[RiAdUKF; Algorithm 21 of \cite{Menegaz2016}]
\label{alg:RiAdUKF}Set the initial estimates $\est{\man{\state}}_{0|0}:=\mean{\man{\state}}_{0}$
and $\est{\man{\cov}}_{\state\state}^{0|0}:=\man{\cov}_{\state\state}^{0}$.
For $k=1,...,k_{f}$, perform the following steps:
\begin{enumerate}[labelsep=0.1cm,leftmargin=0.45cm]
\item \label{enu:RiAdUKF-State-prediction}State prediction. \textup{
\begin{align}
\Big(\est{\man{\state}}_{k|k-1}^{*},\est{\man{\cov}}_{\man{\state\state},*}^{k|k-1}\Big) & :=\riUT_{1}\Big(f_{k},\est{\man{\state}}_{k-1|k-1},\est{\man{\cov}}_{\man{\state\state}}^{k-1|k-1}\Big)\label{eq:RiAdUKF-RiUT1}\\
\est{\man{\state}}_{k|k-1} & :=\rexp_{\est{\man{\state}}_{k|k-1}^{*}}\mean{\pnoise}_{k}\label{eq:RiAdUKF-MeanProcNoise}\\
\est{\man{\cov}}_{\man{\state\state}}^{k|k-1} & :=\est{\man{\cov}}_{\man{\state\state},*}^{k|k-1}+\covpnoise_{k}.\label{eq:RiAdUKF-CovProcNoise}
\end{align}
}
\item \label{enu:RiAdUKF-Measurement-prediction}Measurement prediction.
\textup{
\begin{multline}
\Big(\est{\man{\meas}}_{k|k-1}^{*},\est{\man{\cov}}_{\man{yy},*}^{k|k-1},\est{\man{\cov}}_{\man{\state y}}^{k|k-1}\Big):=\\
\riUT_{2}\left(h_{k},\est{\man{\state}}_{k|k-1},\est{\man{\cov}}_{\man{\state\state}}^{k|k-1}\right)\label{eq:RiAdUKF-RiUT2}
\end{multline}
\begin{align}
\est{\man{\meas}}_{k|k-1} & :=\rexp_{\est{\man{\meas}}_{k|k-1}^{*}}\mean{\mnoise}_{k}\label{eq:RiAdUKF-MeanMeasNoise}\\
\est{\man{\cov}}_{\man{yy}}^{k|k-1} & :=\est{\man{\cov}}_{\man{yy},*}^{k|k-1}+\covmnoise_{k}.\label{eq:RiAdUKF-CovMeasNoise}
\end{align}
}
\item \label{enu:RiAdUKF-State-correction}State correction.\textup{
\begin{align}
\man{\kgain}_{k} & :=\est{\man{\cov}}_{\man{\state}\man{\meas}}^{k|k-1}\Big(\est{\man{\cov}}_{\man{\meas}\man{\meas}}^{k|k-1}\Big)^{-1}\label{eq:RiAdUKF-KalmanGain}\\
\hat{\state}_{k|k}^{TM} & :=\man{\kgain}_{k}\rlog_{\est{\man{\meas}}_{k|k-1}}\big(\outcome{\man{\meas}}_{k}\big)\label{eq:RiAdUKF-tangent-correct-estimate}\\
\est{\man{\state}}_{k|k} & :=\rexp_{\est{\man{\state}}_{k|k-1}}\big(\hat{\state}_{k|k}^{TM}\big)\label{eq:RiAdUKF-corrected-estimate}\\
\est{\man{\cov}}_{\man{\state}\man{\state}}^{k|k,\est{\man{\state}}_{k|k-1}} & :=\est{\man{\cov}}_{\man{\state}\man{\state}}^{k|k-1}-\man{\kgain}_{k}\est{\man{\cov}}_{\man{\meas}\man{\meas}}^{k|k-1}\man{\kgain}_{k}^{T}\nonumber \\
\est{\man{\cov}}_{\man{\state}\man{\state}}^{k|k} & :=\partransp\Big(\est{\man{\cov}}_{\man{\state}\man{\state}}^{k|k,\est{\man{\state}}_{k|k-1}},\est{\man{\state}}_{k|k-1},\est{\man{\state}}_{k|k}\Big).\label{eq:RiAdUKF-corrected-covariance}
\end{align}
}
\end{enumerate}
\end{lyxalgorithm}
A\emph{ll steps of the} \emph{RiUKFs are justified by and coherent
with the other results of this work}. Among these, the most important
are $\sr$, $\riUT$ and Riemannian systems.

The notations $\riUT_{1}$ and $\riUT_{2}$ {[}in (\ref{RiAuUKF-RiUT1}),
(\ref{RiAuUKF-RiUT2}), (\ref{eq:RiAdUKF-RiUT1}), and (\ref{eq:RiAdUKF-RiUT2}){]}
indicate these $\riUT$s \emph{can have different forms}. The output
of $\riUT_{1}$ has only two terms—which is different from the number
of mapped variables in Definition \ref{def:Riemannian-Unscented-Transformation}—meaning
that only the first two variables of the output of Definition \ref{def:Riemannian-Unscented-Transformation}
are needed.

We can consider \emph{not regenerating the independent set of $\riUT_{2}$}
when $\riUT_{1}=\riUT_{2}$. Let $\man{\spa}_{*}^{k|k-1}$ be the
\emph{dependent set of $\riUT_{1}$} and $\man{\spa}_{i}^{k|k-1}$
the \emph{dependent set of $\riUT_{2}$}. Because, from (\ref{RiAuUKF-RiUT2})
and (\ref{eq:RiAdUKF-RiUT2}), $\man{\spa}_{i,*}^{k|k-1}$ and $\man{\spa}^{k|k-1}$
are different objects, we say $\man{\spa}_{*}^{k|k-1}$ is regenerated.
Nonetheless, we could set $\man{\spa}_{i,*}^{k|k-1}=\man{\spa}^{k|k-1}$;
consequently, the computational effort of the filter would decrease—calculating
a new $\man{\spa}^{k|k-1}$ can be computationally because it includes
calculating a square-root matrix of $\est{\man{\cov}}_{\man{\state\state},aug}^{k|k-1}$
or $\est{\man{\cov}}_{\man{\state\state}}^{k|k-1}$. But in this case,
i) the estimation quality of the RiAdUKF would possibly deteriorate—it
has been shown for the Euclidean case (cf. Section 5.1 of \cite{Menegaz2016})—and
ii) the reasoning behind the RiUKFs explained in the second paragraph
of this section would not be true anymore.

After choosing the manifolds' atlases, \emph{all} expressions for
the Riemannian exponentials, logarithms, etc., \emph{must be} \emph{coherent}
with the chosen parameterizations. These transformations, as well
as other elements in these filters such as covariances, have different
expressions \emph{depending on the parameterizations} defining the
manifolds.

We can find $\risr$s (with $\weight_{i}^{m}>0$ for every $i=1,\ldots,\nsp$)
\emph{by first finding $\sr$s in tangent spaces} (see the last paragraph
of Section \ref{sec:Riemannian-sigma-representations}). The independent
sets of $\riUT_{1}$ and $\riUT_{2}$ can be difficult to find. Fortunately,
closed forms of $\risr$s (which can be independent sets of $\riUT s$)
can be found from closed forms of normalized $\sr$s by using Theorem
\ref{thm:Euclidean-to-Riemannian-sigma-rep}. 

The method for obtaining the sample means of $\riUT_{1}$ and $\riUT_{2}$
affects the computation efforts of the RiUKFs because, following \cite{Pennec2006},
we define these sample means as optimization problems (Section \ref{sec:Statistics-of-weighted}).
Sometimes there exist closed forms, but more often it requires optimization
algorithms. The reader will find efficient options in \cite{Absil2008,Pennec2006,Moakher2002,Pennec1998}
and in the MATLAB and Python toolbox ManOpt \cite{Boumal2014}\footnote{Available for download at \textcolor{blue}{\url{https://www.manopt.org/}}.}.

Computational efforts of the RiUKFs also varies with the underlying
manifolds and their atlases because the expressions for exponentials,
logarithms and parallel transports change with them. The reader can
also refer to the ManOpt toolbox for many efficient implementations
of these operations.

Apart from these three factors, computational efforts majorly depends
on the square-rooting involved in the $\text{Ri}\sr$ calculations
and the $\est{\man{\cov}}_{\man{\meas}\man{\meas}}^{k|k-1}$ inversion
in the Kalman gain calculations. Since we can find $\risr$s\textit{
}\textit{\emph{by finding $\sr$s in tangent spaces and, to the best
of our knowledge, all known $\sr$s require square-rooting a covariance
matrix (cf. \cite{Menegaz2015}), the computational complexity of
these operations in }}(\ref{RiAuUKF-RiUT1}) is $\mathcal{O}([\dimstate+\dimpnoise]^{3})$,
in (\ref{RiAuUKF-RiUT2}) $\mathcal{O}([\dimmeas+\dimmnoise]^{3})$,\textit{\emph{
in }}(\ref{eq:RiAdUKF-RiUT1}) is $\mathcal{O}(\dimstate^{3})$, and
in (\ref{eq:RiAdUKF-RiUT2}) $\mathcal{O}(\dimmeas^{3})$. The \textit{\emph{computational
complexity}} of the $\est{\man{\cov}}_{\man{\meas}\man{\meas}}^{k|k-1}$
inversion is $\mathcal{O}(\dimmeas^{3})$ in both (\ref{eq:RiAuUKF-KalmanGain})
and (\ref{eq:RiAdUKF-KalmanGain}).

RiUKFs are\emph{ generalizations of UKFs}. Every UKF is a RiUKF, and
every RiUKF for Euclidean state-variables is a UKF. It is easy to
see that, if $\mana_{\man{\state}}$ and $\mana_{\man{\meas}}$ are
Euclidean spaces, then RiAuUKF is equivalent to AuUKF (Algorithm 7
of \cite{Menegaz2016}), and RiAdUKF to AdUKF (Algorithm 6 of \cite{Menegaz2016}).

Since Cartesian products of Riemannian manifolds are also Riemannian
manifolds (e.g., $\sset 3\times\realset^{n}$) \cite{DoCarmo1992},
the proposed RiUKF also estimates systems with state variables belonging
to \emph{Cartesian products} of Riemannian manifolds.

The Kalman gain $\man{\kgain}_{k}$ in (\ref{eq:RiAuUKF-KalmanGain})
and (\ref{eq:RiAdUKF-KalmanGain}) could be defined in a more general
way, as done in (\ref{eq:UKF-Riemannian-CartesianSolution-3}). However,
it would imply more computational effort—the dimension of the sigma
points and matrices would be higher—at the exchange of no advantage,
at least at present; perhaps benefits can be obtained from (\ref{eq:UKF-Riemannian-CartesianSolution-3})
in future works.

The three assumptions cited at the beginning of Section \ref{sec:Riemannian-sigma-representations}
impose some limitations on the RiUKFs. Assumption \ref{enu:assumption1}
limits the RiUKFs to the case of geodesically-complete Riemannian
manifolds: still there are many of these manifolds useful for practical
applications, such as unit spheres, special orthogonal groups, special
Euclidean groups, real projective spaces, special unitary groups,
Grassmann manifolds, among others (cf. \cite{Absil2008} and Section
\ref{subsec:Kalman-filtering-in}). Assumption \ref{enu:assumption2}
imposes careful choice of $\man{\cov}_{\man{\state}\man{\state}}^{0}$,
$\man{\covpnoise}_{k}$, $\man{\covmnoise}_{k}$ (or $Q_{k}$ and
$\covmnoise_{k}$ for the RiAdUKF): their values should be consistent
with the logarithms in their definitions {[}or in (\ref{eq:Riemannian-additive-system})
in the case of the RiAdUKF{]}; since these covariances are tuning
parameters and are often set based on intuition, an user could chose
inconsistent (too great) values; this would probably result on either
inconsistent sigma points—because the tangent sigma points would be
outside the tangent cut locus—or on some divergence in the algorithm,
such as non-positive state covariance matrix. Assumption \ref{enu:assumption6}
will not, in most cases, impose other limitations if the user model
the system equations and parameters consistently.

We can find \emph{particular cases} of RiUKFs by choosing particular
forms of $\risr$s; Table \ref{table:Riemannian-augmented-minimum-kalman-filters}
shows some cases for $\riUT_{1}=\riUT_{2}$—the second and third columns
contain the filters. Each filter is the resulting variant of using
i) the corresponding RiUKF in the \emph{heading row} of its column
(RiAuUKF or RiAdUKF), and ii) the corresponding $\risr$ written in
the first column of its row. For instance, the Riemannian Minimum
AuUKF (RiMiAuUKF in the first row and second column), is the result
of the RiAuUKF with the RiMi$\sigma$R (Corollary \ref{cor:Riemannian-sr-minimum-numbers}).
All filters in Table \ref{table:Riemannian-augmented-minimum-kalman-filters}
are \emph{new}.

\begin{table}
\caption{RiUKF Variants for some Ri$\sigma$Rs.\label{table:Riemannian-augmented-minimum-kalman-filters}}
\begin{centering}
\begin{tabular}{ccc}
\textbf{$\sigma$R}\footnotemark[1] & \textbf{AuUKF}\footnotemark[1] & \textbf{AdUKF }\footnotemark[1]\tabularnewline
\hline 
\textbf{RiMi$\sigma$R} & RiMiAuUKF & RiMiAdUKF\tabularnewline
\textbf{RiRhoMi$\sigma$R} & RiRhoMiAuUKF & RiRhoMiAdUKF\tabularnewline
\textbf{RiMiSy$\sigma$R} & RiMiSyAuUKF & RiMiSyAdUKF\tabularnewline
\textbf{RiHoMiSy$\sigma$R} & RiHoMiSyAuUKF & RiHoMiSyAdUKF\tabularnewline
\end{tabular}
\par\end{centering}
\footnotemark[1]{Ad for Additive, Au for Augmented, Ho for Homogeneous,
Mi for Minimum, Ri stands for Riemannian, $\sigma$R for \textbf{$\sigma$-}Representation\textbf{,}
Sy for Symmetric, UKF for Unscented Kalman Filter. Rho stand for Rho
itself; see also the acronyms list in Appendix \ref{subsec:Notation-and-Acronyms}}
\end{table}

An RiUKF for the partially-additive system (\ref{eq:Riemannian-partiallly-additive-system1})
is given by step \ref{enu:RiAuUKF-State-prediction} of the RiAuUKF
with steps \ref{enu:RiAdUKF-Measurement-prediction} and \ref{enu:RiAdUKF-State-correction}
of the RiAdUKF, and for (\ref{eq:Riemannian-partiallly-additive-system2})
is given by step \ref{enu:RiAdUKF-State-prediction} of the RiAdUKF
with steps \ref{enu:RiAuUKF-Measurement-prediction} and \ref{enu:RiAuUKF-State-correction}
of the RiAuUKF.

For (\ref{eq:Riemannian-additive-system}), (\ref{eq:Riemannian-partiallly-additive-system1})
and (\ref{eq:Riemannian-partiallly-additive-system2}) when either
$f_{k}$ or $h_{k}$ are the identity function, we can simplify their
filters by skipping sigma points calculations; hence saving computation
effort. If, for example, $f_{k}(\man{\state})=\man{\state}$, then
the following two equations can replace the state prediction (e.g.,
the step \ref{enu:RiAdUKF-State-prediction} of the RiAdUKF):
\begin{align*}
\est{\man{\state}}_{k|k-1} & :=\rexp_{\est{\man{\state}}_{k-1|k-1}}\mean{\pnoise}_{k}\\
\est{\man{\cov}}_{\man{\state\state}}^{k|k-1} & :=\est{\man{\cov}}_{\man{\state\state}}^{k-1|k-1}+\covpnoise_{k}.
\end{align*}
The case $\mfunc_{k}(\man{\state})=\man{\state}$ is similar.

\subsection{Relation with the literature}

\label{subsec:Relation-with-the}

To the best of our knowledge, the \emph{UKF for Riemannian manifolds}
(\emph{UKFRM})\textit{ }\textit{\emph{of}} \textit{\emph{\cite{Hauberg2013}}}
is the only UKF for any geodesically-complete Riemannian manifold
in the literature. Consider system (\ref{eq:Riemannian-general-system})
and define the following functions—cf. (1) and (2) of \cite{Hauberg2013}—:
\begin{equation}
\pfunc_{k}^{*}(\man{\state}_{k-1}):=\pfunc_{k}(\man{\state}_{k-1},\man{\pnoise}_{k-1}),\,\mfunc_{k}^{*}(\man{\state}_{k}):=\mfunc_{k}(\man{\state}_{k},\man{\mnoise}_{k}).\label{eq:system-of-hauberg}
\end{equation}
Suppose that i) the initial state $\man{\state}_{0}$ is characterized
by $\man{\state}_{0}\sim(\mean{\man{\state}}_{0},\man{\cov}_{\man{\state}\man{\state}}^{0})_{\mana_{\man{\state}}},$
and ii) the measurements $\outcome{\man{\meas}}_{1}$, $\outcome{\man{\meas}}_{2}$,
..., $\outcome{\man{\meas}}_{k_{f}}$ are given. Let 
\[
\mbox{HoMiSy}\sr:(\mean{\rva},\cov_{\rva\rva})\mapsto\{\spa_{i,},\weight_{i}\}_{i=1}^{\nsp}
\]
be a function mapping the mean $\mean{\rva}$ and covariance $\cov_{\rva\rva}$
of a given random vector $\rva$ to a HoMiSy$\sigma$R (Corollary
3 of \cite{Menegaz2015}). Then the UKFRM of \cite{Hauberg2013} is
given by the following algorithm:
\begin{lyxalgorithm}[UKFRM of \cite{Hauberg2013}]
\label{alg:UKF-Hauberg}Set $\nsp:=2\dimstate+1$ and the initial
estimates $\est{\man{\state}}_{0|0}:=\mean{\man{\state}}_{0}$ and
$\est{\man{\cov}}_{\state\state}^{0|0}:=\man{\cov}_{\state\state}^{0}$.
For $k=1,...,k_{f}$, perform the following steps:
\end{lyxalgorithm}
\begin{enumerate}[labelsep=0.1cm,leftmargin=0.45cm]
\item State prediction. 
\begin{align}
 & \big\{\spa_{i,k-1|k-1}^{TM},\weight_{i}\big\}{}_{i=1}^{\nsp}:=\mbox{HoMiSy}\sr\Big([0]_{\dimstate},\est{\man{\cov}}_{\man{\state}\man{\state}}^{k-1|k-1}\Big)\label{eq:UKF-Hauberg-tangent-previous-sigma-rep}\\
 & \man{\spa}_{i}^{k-1|k-1}:=\exp_{\est{\man{\state}}_{k-1|k-1}}\big(\spa_{i,k-1|k-1}^{TM}\big),\,i=1,\ldots,\nsp\label{eq:UKF-Hauberg-previous-sigma-rep}\\
 & \man{\spa}_{i,*}^{k|k-1}:=\pfunc_{k}^{*}\big(\man{\spa}_{i}^{k-1|k-1}\big),\,i=1,\ldots,\nsp\nonumber \\
 & \est{\man{\state}}_{k|k-1}:=\arg\min_{\man a\in\mana_{\man{\state}}}\sum_{i=1}^{\nsp}\weight_{i}\dist^{2}\big(\man{\spa}_{i,*}^{k|k-1},\man a\big)\nonumber \\
 & \est{\man{\cov}}_{\man{\state}\man{\state}}^{k|k-1}:=\sum_{i=1}^{\nsp}\weight_{i}\Big(\rlog_{\est{\man{\state}}_{k|k-1}}\big(\man{\spa}_{i,*}^{k|k-1}\big)\Big)\Big(\diamond\Big)^{T}.\label{eq:UKF-Hauberg-Predicted-X-COV}
\end{align}
\item Measurement prediction.
\begin{align}
 & \Big\{\spa_{i,k|k-1}^{TM},\weight_{i}\Big\}_{i=1}^{\nsp}:=f_{k}\mbox{HoMiSy}\sr\Big([0]_{\dimstate},\est{\man{\cov}}_{\man{\state}\man{\state}}^{k|k-1}\Big)\label{eq:UKF-Hauberg-tangent-predicted-sigma-rep}\\
 & \man{\spa}_{i}^{k|k-1}:=\exp_{\est{\man{\state}}_{k|k-1}}\big(\spa_{i,k|k-1}^{TM}\big),\,i=1,\ldots,\nsp\label{eq:UKF-Hauberg-second-predicted-sigma-rep}\\
 & \man{\spb}_{i}^{k|k-1}:=h_{k}^{*}\big(\man{\spa}_{i}^{k|k-1}\big),\,i=1,\ldots,\nsp\nonumber \\
 & \est{\man{\meas}}_{k|k-1}:=\arg\min_{\man b\in\mana_{\man{\meas}}}\sum_{i=1}^{\nsp}\weight_{i}\dist^{2}\left(\man{\spb}_{i}^{k|k-1},\man b\right)\nonumber \\
 & \est{\man{\cov}}_{\man{\meas}\man{\meas}}^{k|k-1}:=\sum_{i=1}^{\nsp}\weight_{i}\Big(\rlog_{\est{\man{\meas}}_{k|k-1}}\big(\man{\spb}_{i}^{k|k-1}\big)\Big)\Big(\diamond\Big)^{T}\label{eq:UKF-Hauberg-Predicted-Y-COV}\\
 & \est{\man{\cov}}_{\man{\state}\man{\meas}}^{k|k-1}:=\sum_{i=1}^{\nsp}\weight_{i}\Big(\rlog_{\est{\man{\state}}_{k|k-1}}\big(\man{\spa}_{i}^{k|k-1}\big)\Big)\nonumber \\
 & \quad\quad\quad\quad\quad\quad\quad\quad\quad\quad\quad\quad\quad\Big(\rlog_{\est{\man{\meas}}_{k|k-1}}\big(\man{\spb}_{i}^{k|k-1}\big)\Big)^{T}.\nonumber 
\end{align}
\item State correction.
\begin{align}
\man{\kgain}_{k} & :=\est{\man{\cov}}_{\man{\state}\man{\meas}}^{k|k-1}\Big(\est{\man{\cov}}_{\man{\meas}\man{\meas}}^{k|k-1}\Big)^{-1}\nonumber \\
\hat{\state}_{k|k}^{TM} & :=\hat{\state}_{k|k-1}^{TM}+\man{\kgain}\rlog_{\est{\man{\meas}}_{k|k-1}}\big(\outcome{\man{\meas}}_{k}\big)\label{eq:ukf-hauberg-state-corrected-estimate}\\
\est{\man{\state}}_{k|k} & :=\rexp_{\est{\man{\state}}_{k|k-1}}\big(\hat{\state}_{k|k}^{TM}\big)\nonumber \\
\est{\man{\cov}}_{\man{\state}\man{\state}}^{k|k,\est{\man{\state}}_{k|k-1}} & :=\est{\man{\cov}}_{\man{\state}\man{\state}}^{k|k-1}-\man{\kgain}_{k}\est{\man{\cov}}_{\man{\meas}\man{\meas}}^{k|k-1}\man{\kgain}_{k}^{T}\nonumber \\
\est{\man{\cov}}_{\man{\state}\man{\state}}^{k|k} & :=\partransp\Big(\est{\man{\cov}}_{\man{\state}\man{\state}}^{k|k,\est{\man{\state}}_{k|k-1}},\est{\man{\state}}_{k|k-1},\est{\man{\state}}_{k|k}\Big).\nonumber 
\end{align}
\end{enumerate}

Compared with the UKFRM of \cite{Hauberg2013}, we can point out the
following five improvements of the RiUKFs:

\begin{enumerate}[labelsep=0.1cm,leftmargin=0.45cm]
\item The noises are incorporated into the RiUKFs, but in the UKFRM they
are not. In the RiAuUKF, the noises are incorporated by realizing
the augmented sigma points in the process and measurement functions
{[}equations (\ref{RiAuUKF-RiUT1}) and (\ref{RiAuUKF-RiUT2}){]};
and in the RiAdUKF, by ``adding'' (in the tangent space) their means
and covariances {[}equations (\ref{eq:RiAdUKF-MeanProcNoise}), (\ref{eq:RiAdUKF-CovProcNoise}),
(\ref{eq:RiAdUKF-MeanMeasNoise}), (\ref{eq:RiAdUKF-CovMeasNoise}){]}.\\
However, the UKFRM exclude the noises. Even though the UKFRM of \cite{Hauberg2013}
considers a system with process and measurement noises {[}cf. (\ref{eq:system-of-hauberg}){]},
they do \emph{not} influence any estimate within the UKFRM; these
noises' statistics \emph{do not} \emph{appear at any step }of the
UKFRM—commonly, filters consider these statistics when calculating
the predicted covariances, but this is also not the case for the UKFRM
{[}cf. (\ref{eq:UKF-Hauberg-Predicted-X-COV}) and (\ref{eq:UKF-Hauberg-Predicted-Y-COV}){]}.\\
We can point out at least two consequences of this absence of the
noise elements:
\begin{enumerate}
\item the Euclidean case of the UKFRM is \emph{not equivalent} \emph{to
any} (Euclidean) UKF. This can be seen by considering Euclidean manifolds
in Algorithm \ref{alg:UKF-Hauberg} (cf. the last paragraph of Sections
\ref{subsec:Statistics-of-random} and of \ref{sec:Statistics-of-weighted}).
Besides, to the best of our knowledge, there is no UKF without process
and measurement noises covariance (cf. \cite{Menegaz2015,Menegaz2016}). 
\item the UKFRM might \emph{diverge} in situations in which the RiUKFs do
not. This behavior can be seen in the following simple example: consider
(\ref{eq:Riemannian-additive-system}) and (\ref{eq:system-of-hauberg})
with $\mana_{\man{\state}}^{\dimstate}=\mana_{\man{\state}}^{\dimpnoise}=\mana_{\man{\state}}^{\dimmeas}=\mana_{\man{\state}}^{\dimmnoise}=\realset$.
Suppose that i) the initial state is $\state_{0}\sim(1,1)_{\realset},$
ii) the noise covariances are $\covpnoise_{k}=\covmnoise_{k}=1$,
iii) the system functions are $\pfunc_{k}(\state_{k-1})=\pfunc_{k}^{*}(\state_{k-1})=\state_{k-1}$
and $\mfunc_{k}(\state_{k})=\mfunc_{k}^{*}(\state_{k})=1-\state_{k}$,
and iv) the measurements are $\outcome{\meas}_{1}=\cdots=\outcome{\meas}_{k_{f}}=1$.
For this example, we ran the (linear) KF (cf. \cite{Jazwinsky1970}),
the RiAdUKF, and the UKFRM. Both the KF and the RiAdUKF provided the
same estimates, but the UKFRM did not provide consistent results;
the simulation was halted because the corrected covariance ($\est{\man{\cov}}_{\man{\state}\man{\state}}^{2|2}$)
lost its positiveness. Similar results occurred in the simulations
of Section \ref{sec:Numerical-Example}.
\end{enumerate}
\item We introduced a consistent definition {[}equation (\ref{eq:Riemannian-additive-system}){]}
for the system associated with the RiAdUKF. To the best of our knowledge,
(\ref{eq:Riemannian-additive-system}) is the first consistent additive-noise
Riemannian stochastic discrete-time dynamic system.
\item To the best of our knowledge, the RiUKFs are the \emph{first} UKFs
for Riemannian state-space systems considering noises with \emph{non-zero
means}. Even for simple manifolds such as the unit sphere, we could
not find a UKF considering this case.
\item \emph{A}ll the equations of our RiUKFs are \emph{formally justified}.
These justifications are the following ones:
\begin{enumerate}
\item \emph{The equations of steps }\ref{enu:RiAuUKF-State-prediction}
and \ref{enu:RiAuUKF-Measurement-prediction} of the RiUKFs\emph{
are justified by Definition }\ref{def:Consider-a-Riemannian},\emph{
Theorem }\ref{thm:Euclidean-to-Riemannian-sigma-rep} and Corollary
\ref{sec:Riemannian-Unscented-Transformat}.
\item \emph{Equations }(\ref{eq:RiAuUKF-KalmanGain})\emph{ and }(\ref{eq:RiAdUKF-KalmanGain})\emph{
(the Kalman Gains) are justified in Section \ref{subsec:State-and-measurement=00005B}}.
This form of the Kalman gain $\man{\kgain}_{k}$ in (\ref{eq:RiAuUKF-KalmanGain})
and (\ref{eq:RiAdUKF-KalmanGain}) follows as a particular case of
the Kalman gain of a more general system ($\man{\kgain}_{k,**}$)
where the state and the measurement belong to the product $\mana_{\man{\state}}\times\mana_{\man y}$.
\item \emph{The equations of step }\ref{enu:RiAuUKF-State-correction} of
the RiUKFs\emph{ are justified in Section \ref{subsec:Correction-equations}}.
We showed that they follow from considering i) $\state_{k|k-1}^{TM}$
and $\meas_{k|k-1}^{TM}$ normally-joint distributed {[}equation (\ref{eq:predicted-tangential-statistics-are-normal}){]},
and ii) $\state_{k|k}^{TM}$ given by a linear correction of $\state_{k|k-1}^{TM}$
by $(\outcome{\meas}_{k}^{TM}-\meas_{k|k-1}^{TM})$ {[}equation (\ref{eq:correction-equation-of-xk-tangent}){]}. 
\end{enumerate}
\item (Euclidean) UKFs are particular cases of the RiUKFs (cf. Section \ref{subsec:Riemannian-Unscented-Filters}).
\end{enumerate}
Altogether, we can say the RiUKFs have novelties compared with the
UKF for Riemannian state-space systems of the literature.

\section{Example: Satellite Attitude Tracking}

\label{sec:Numerical-Example}

In this section, we apply the developed theory to estimate the attitude
of a satellite in a realistic scenario (cf. \cite{Crassidis2007}).

The set of possible attitudes of a rotating body is not a Euclidean
space, but a three dimensional smooth manifold known as $SO(3)$.
This manifold has many different topological properties from a Euclidean
space: for instance, it is compact whilst Euclidean spaces are not.
Due to this difference, Euclidean UKFs designed over Euclidean spaces
may not work properly: its estimates may not stay within the state-space
manifold, resulting in poor performance and poor accuracy \cite{Crassidis2003}. 

Although we could apply an RiUKFs for $\SO 3$ in this example, we
prefer to apply an RiUKF for the set of unit quaternions $\sset 3$
because they represent, without singularities \cite{Stuelpnagel1964},
attitudes using the minimal set of parameters. Let $\quat q_{i}=\begin{bmatrix}\eta_{i} & \myvec{\epsilon}_{i}^{T}\end{bmatrix}^{T}\in\mathbb{R}^{4}$,
where $\eta_{i}\in\mathbb{R}$ and $\myvec{\epsilon}_{i}\in\mathbb{R}^{3}$.
It is possible to prove that the three dimensional sphere
\begin{equation}
S^{3}=\{(q_{1},q_{2},q_{3},q_{4})\in\mathbb{R}^{4}:q_{1}^{2}+q_{2}^{2}+q_{3}^{2}+q_{4}^{2}=1\}\label{eq:sphere_manifold}
\end{equation}
is a Riemannian manifold and the product 
\[
\quat q_{1}\otimes\quat q_{2}=\begin{bmatrix}\eta_{1}\eta_{2}-\myvec{\epsilon}_{1}^{T}\myvec{\epsilon}_{2}\\
\eta_{1}\myvec{\epsilon}_{2}+\eta_{2}\myvec{\epsilon}_{1}+\myvec{\epsilon}_{1}\times\myvec{\epsilon}_{2}
\end{bmatrix}.
\]
is closed. For a rotation of an angle $\theta$ around an unit vector
$n$, there are two associated unit quaternions $\quat q$ and $\quat q'$
such that
\[
\quat q=\cos\left(\frac{\theta}{2}\right)+\imvec\mathbf{n}\sin\left(\frac{\theta}{2}\right),\quad\quat q'=-\quat q.
\]

Let $\quata(t)\in\sset 3$ be the attitude of the satellite at the
time instant $t$, and $\omega(t)\in\realset^{3}$ its the angular
velocity. The evolution of $\quata(t)$ over time can be described
by the following differential equation \cite{Zipfel2007}: 

\begin{equation}
\dot{\quata}\left(t\right)=\frac{1}{2}\quat{\omega}(t)\otimes\quata\left(t\right),\label{eq:satellite-system}
\end{equation}
where $\quat{\omega}\in\mathbb{R}^{4}$ is given by $\quat{\omega}=\begin{bmatrix}0 & \omega^{T}\end{bmatrix}^{T}$. 

We generate synthetic data by a fourth order Runge-Kutta integration
of (\ref{eq:satellite-system}) over the interval $[0\text{s},20\text{s}]$
with angular velocity
\[
\omega\left(t\right)=\left[\begin{array}{c}
0.03\sin\left(\left[\pi t/600\right]\degree\right)\\
0.03\sin\left(\left[\pi t/600\right]\degree-300\degree\right)\\
0.03\sin\left(\left[\pi t/600\right]\degree-600\degree\right)
\end{array}\right]
\]
and initial state $\quata\left(0\right)=0.96+\imvec[0.13,\,0.19,\,\sqrt{1-0.96^{2}-0.13{}^{2}-0.19{}^{2}}]^{T}$
.

For filtering, we consider (\ref{eq:Riemannian-additive-system})
with $\man{\state}_{k}=\quata(k\samptime)$ 
\begin{align*}
\theta(t) & :=\norm{\omega(t)}\frac{\samptime}{2}\\
\pfunc_{k}\big(\man{\state}_{k-1}\big) & =\begin{bmatrix}\cos\theta(t) & \frac{\omega^{T}(t)}{\norm{\omega(t)}}\sin\theta(t)\end{bmatrix}^{T}\quatmultsymbol\quat{\state}_{k-1}\\
\mfunc_{k}\big(\man{\state}_{k}\big) & =\man{\state}_{k},
\end{align*}
$\mean{\pnoise}_{k}=\mean{\mnoise}_{k}=[0]_{3\times1}$, $\boldsymbol{\covpnoise}_{k}=(0.31236\times10^{-6})^{2}I_{3}$,
and $\boldsymbol{\covmnoise}_{k}=(0.5\pi/180\times10^{-6})^{2}I_{3}$.
These values for $\boldsymbol{\covpnoise}_{k}$ and $\boldsymbol{\covmnoise}_{k}$
were chosen according to \cite{Crassidis2003}.

We performed $1,000$ simulations with the RiUKFs of Table \ref{table:Riemannian-augmented-minimum-kalman-filters}
and the UKFRM of \cite{Hauberg2013}. To calculate Riemannian means,
we used the gradient descent method of \cite{Pennec1998} with a threshold
of $10^{-6}$; and for Riemannian exponentials, Riemannian logarithms,
and parallel transport, we used the MATLAB toolbox ManOpt \cite{Boumal2014}.

For all simulations, the RiUKFs of Table \ref{table:Riemannian-augmented-minimum-kalman-filters}
provided good estimates, with a Root Mean Square Error in the order
of $10^{-6}$ (Table \ref{tab:RMSE-stat-track}). The RiMiAdUKF or
the RiRhoMiAdUKF are the best alternatives for this example because
i) it demands less computational effort than the other filters—it
is additive and is composed of the least number of sigma points (cf.
Corollary \ref{cor:Riemannian-sr-minimum-numbers})— and ii) all RiUKFs
performed almost equally. 

\begin{table}
\centering{}\caption{Root Mean Square Error ($\times10^{-6}$) of each RiUKF in Table \ref{table:Riemannian-augmented-minimum-kalman-filters}
considering 1,000 simulations of a satellite attitude tracking example.\label{tab:RMSE-stat-track}}
\begin{tabular}{cccc}
\hline 
RiMiAuUKF & RiRhoMiAuUKF & RiMiSyAuUKF & RiHoMiSyAuUKF\tabularnewline
2,612 & 2,614 & 2,614 & 2,614\tabularnewline
\hline 
RiMiAdUKF & RiRhoMiAdUKF & RiMiSyAdUKF & RiHoMiSyAdUKF\tabularnewline
2,612 & 2,613 & 2,613 & 2,613\tabularnewline
\hline 
\end{tabular}
\end{table}

The UKFRM failed in all the $1,000$ simulations; in every simulation,
the state covariance estimate lost its positiveness. Nonexistence
of noise terms in the UKFRM might explain this problematic behavior
(cf. Section \ref{subsec:Relation-with-the}).

\section{Conclusions}

\label{sec:Conclusions}

In this work, we extend the systematization of the Unscented Kalman
Filtering theory we developed in \cite{Menegaz2015} towards estimating
the state of Riemannian systems. In this systematization, we introduce
the following results\footnote{These results were first presented in Menegaz's PhD thesis \cite{Menegaz2016}.}
(all results are mathematically justified):
\begin{enumerate}
\item A Riemannian extension of the $\sigma$-representation ($\sr$ ):
the Riemannian $\sigma$-representation ($\risr$, Section \ref{sec:Riemannian-sigma-representations}).
\item A technique to obtain closed forms of the $\risr$ by closed forms
of the $\sr$ (Theorem \ref{thm:Euclidean-to-Riemannian-sigma-rep}).
Using this result, we discover (Corollary \ref{cor:Riemannian-sr-minimum-numbers})
\begin{enumerate}
\item the minimum number of sigma points of an $\risr$,
\item the minimum number of sigma points of a symmetric $\risr$,
\item closed forms for the minimum $\risr$, and
\item closed forms for the minimum symmetric $\risr$. 
\end{enumerate}
\item An \textit{\emph{additive-noise}} Riemannian\textit{ }system definition
(Section \ref{subsec:Riemannian-Systems}). We require this definition
to introduce additive-noise Riemannian UKFs.
\item Kalman correction equations on Riemannian manifolds (Section \ref{subsec:Correction-equations}).
\item New discrete-time Riemannian UKFs (RiUKFs), namely the Riemannian
Additive UKF and the Riemannian Augmented UKF (Section \ref{subsec:Riemannian-Unscented-Filters}).
Besides, we
\begin{enumerate}
\item provide a list of particular variants of these filters (Table \ref{table:Riemannian-augmented-minimum-kalman-filters});
all these variants are new. Compared with the literature's UKF for
Riemannian manifolds (in \cite{Hauberg2013}), our RiUKFs are more
consistent, formally-principled, and general.
\item numerically compare all these particular variants with the literature's
UKF on Riemannian manifolds in a satellite attitude tracking scenario.
For all 1,000 simulations, the new variants provided good estimates,
but the literature's filter diverged; in every simulation, the state
covariance estimate lost its positiveness.
\end{enumerate}
\end{enumerate}

With this work, we hope to have expanded the literature's knowledge
on Kalman filtering and provided a tool for the research community
to improve the performance and stability of many UKFs.

Following this study, we recommend the research community searching
for computationally-implementable variants of RiUKFs. Since concepts
of the Riemannian manifold theory can be very abstract, depending
on the underlying manifold, developing RiUKFs variants is not trivial.

This task is even harder without a generalizing base theory: that
is one of the reasons why, in this work, we develop a general consistent
systematized theory of Unscented Kalman Filters for Riemannian State-Space
Systems. 


\appendix

\section{Appendix}

\subsection{Results relative to Riemannian manifolds \label{appendix:Riemannian-manifolds}}

In this appendix, we provide some results relative to the theory of
Riemannian manifolds. These definitions are mainly based on \cite{DoCarmo1992}.

\begin{defn}[Differentiable manifold \cite{DoCarmo1992}]
\label{def:differentiable-manifold}A \emph{differentiable manifold}
of dimension $n$ is a pair $\left(\mana,\atlasa\right)$ where $\mana$
is a set, and $\atlasa=\{(U_{\vecconsta},\charta_{a})\}$, called
atlas, a family of injective mappings (\emph{charts}) $\charta_{\vecconsta}:U_{\vecconsta}\subset\realset^{n}\rightarrow\mana$
of open sets $U_{\vecconsta}$ of $\realset^{n}$ into $\mana$ such
that:
\begin{enumerate}
\item \label{enu:def_dif_man_item1}$\bigcup_{a}\charta_{\vecconsta}(U_{\vecconsta})=\mana.$
\item \label{enu:def_dif_man_item2}for any pair $\vecconsta,\vecconstb$,
with $\charta_{\vecconsta}(U_{\vecconsta})\cap\charta_{\vecconstb}(U_{\vecconstb})=:W\neq\emptyset,$
the sets $\charta_{\vecconsta}^{-1}(W)$ and $\charta_{\vecconstb}^{-1}(W)$
are open sets in $\realset^{n}$, and the mappings $\charta_{\vecconstb}^{-1}\circ\charta_{\vecconsta}$
and $\charta_{\vecconstb}^{-1}\circ\charta_{\vecconsta}$ are differentiable.
\item \label{enu:def_dif_man_item3}The family $\atlasa=\{(U_{\vecconsta},\charta_{\vecconsta})\}$
is \emph{maximal} relative to the conditions \ref{enu:def_dif_man_item1})
and \ref{enu:def_dif_man_item2}).
\end{enumerate}
A pair $(U_{\vecconsta},\charta_{\vecconsta})$ (or the mapping $\charta_{\vecconsta}$)
with $\manpta\in\charta_{\vecconsta}(U_{\vecconsta})$ is called a
\emph{parameterization} of $\mana$ at $\manpta$. For simplicity,
we can denote a differentiable manifold $\left(\mana,\atlasa\right)$
of dimension $n$ by $\mana$ or $\mana^{n}$.
\end{defn}

\begin{defn}[Differentiable function\cite{DoCarmo1992}]
\label{def:differentiable-function}Let $\mana_{1}^{n}$ and $\mana_{2}^{m}$
be differentiable manifolds. A mapping $\manfunca:\mana_{1}\rightarrow\mana_{2}$
is \emph{differentiable at} $\manpta\in\mana_{1}$ if, given a parameterization
$\charta_{2}:V\subset\realset^{m}\rightarrow\mana_{2}$ at $\manfunca(\manpta)$,
there exists a parameterization $\charta_{1}:U\subset\realset^{n}\rightarrow\mana_{1}$
at $\manpta$ such that $\manfunca(\charta_{1}(U))\subset\charta_{2}(V)$
and the mapping
\begin{equation}
\tilde{f}:=\charta_{2}^{-1}\circ\manfunca\circ\charta_{1}:U\subset\realset^{n}\rightarrow\realset^{m}\label{eq:expressionphi-1}
\end{equation}
is differentiable at $\charta_{1}^{-1}(\manpta)$. We say $\manfunca$
is differentiable on an open set of $\mana_{1}$ if it is differentiable
at all of the points of this open set. 
\end{defn}

In this work, \emph{we suppose that all functions are differentiable
unless otherwise stated}.

\begin{defn}[Tangent space \cite{DoCarmo1992}]
\label{def:tangent-space}Let $\mana$ be a differentiable manifold.
A differentiable function $\curve:\interval\rightarrow\mana$ is called
a (differentiable) \emph{curve }in $\mana$. Suppose $\curve(0)=\manpta\in\mana$,
and let $\diffset_{\manpta}(\mana)$ be the set of all functions $\funca:\mana\rightarrow\realset$
that are differentiable at $\manpta$. The \emph{tangent vector to
the curve $\curve$ }at $t=0$ is a function $\curve'(0):\diffset_{\manpta}(\mana)\rightarrow\realset$
given by
\[
\curve'(0)\funca=\left.d(\funca\circ\curve)/dt\right|_{t=0},\quad\funca\in\diffset_{\manpta}(\mana).
\]
Note that $\curve'(0)$ is an operator taking $\funca\in\diffset_{\manpta}(\mana)$
to a scalar $\left.d(\funca\circ\curve)/dt\right|_{t=0}$. A \emph{tangent
vector at $\manpta$ }is a tangent vector of some curve $\curve:\interval\rightarrow\mana$
with $\curve(0)=\manpta$ at $t=0$. The set of all tangent vectors
to $\mana$ at $\manpta$ will be indicated by $T_{\manpta}\mana$.
\end{defn}
The set $T_{\manpta}\mana$ forms a vector space of dimension $n$
and is called the \emph{tangent space of $\mana$ at $\manpta$.}

\begin{defn}[Arc length \cite{Pennec2006}]
\label{def:arc-length}Given an open interval $\interval\subset\realset$,
a differentiable function (Definition \ref{def:differentiable-function})
$\curve:\interval\rightarrow\mana$ is called a (differentiable) \emph{curve
}in $\mana$. Given a curve $\curve$ on $\mana$, the arc length
of $\curve$ in the interval $\left[a,b\right]\subset\interval$ is
defined by
\[
\arclength_{a}^{b}(\curve):=\int_{a}^{b}\norm{\geod'(t)}_{\curve(t)}dt.
\]
\end{defn}

\begin{defn}[Differential of a function]
Let $\mana_{1}$ and $\mana_{2}$ be differentiable manifolds and
$\funca:\mana_{1}\rightarrow\mana_{2}$ a differentiable mapping.
For every $\manpta\in\mana_{1}$ and for each $\tanpta\in T_{\manpta}\mana_{1}$,
choose a differentiable curve $\curve:\interval\rightarrow\mana_{1}$
with $\curve(0)=\manpta,$ $\curve'(0)=\tanpta$. Take $\beta=\funca\circ\curve$.
Then it can be shown that the operator $d\funca_{\manpta}(\tanpta)$
defined by 
\[
d\funca_{\manpta}(\tanpta):=\beta'(0)
\]
is a tangent vector of $T_{\funca(\manpta)}\mana_{2}$. Moreover the
mapping the
\[
df_{\manpta}:T_{\manpta}\mana_{1}\rightarrow T_{\funca(\manpta)}\mana_{2}:\tanpta\mapsto\beta'(0)
\]
is linear and does not depend on the choice of $\curve$ \cite{DoCarmo1992}.
This linear mapping $d\funca_{\manpta}$ is called the \emph{differential}
of $\funca$ at $\manpta$.
\end{defn}

\begin{defn}[Vector field \cite{DoCarmo1992,Absil2008}]
\label{def:vector-field}A \emph{vector field} $\vecfielda$ on a
differentiable manifold $\mana$ is a correspondence that associates
to each point $\manpta\in\mana$ a vector $\vecfielda(\manpta)\in T_{\manpta}\mana$.
Given a vector field $\vecfielda$ on $\mana$ and a differentiable
real-valued function $f:\mana\rightarrow\realset$, we let $\vecfielda f$
denote the real-valued function on $\mana$ defined by 
\begin{eqnarray*}
\left(\vecfielda f\right): & \mana & \rightarrow\realset\\
 & \manpta & \mapsto\tanpta f,\quad\tanpta\in T_{\manpta}\mana.
\end{eqnarray*}
 The set of all vector fields of $\mana$ is denote by $\vfieldset{\mana}$.

The multiplication of a vector field $\vecfielda$ by a function $f:\mana\rightarrow\realset$
is defined by $f\vecfielda$: $\mana\rightarrow T_{\manpta}\mana$:
$\manpta\mapsto f(\manpta)\tanpta$, $\tanpta\in T_{\manpta}\mana$;
and the addition of two vector fields $\vecfielda$ and $\vecfieldb$
by $\vecfielda+\vecfieldb:$ $\mana\rightarrow T_{\manpta}\mana:$
$\manpta\mapsto\vecfielda(\manpta)+\vecfieldb(\manpta)$. The \emph{Lie
bracket of vector fields} is defined as the unique vector field $[\vecfielda,\vecfieldb]$
satisfying $([\vecfielda,\vecfieldb]f)\coloneqq(\vecfielda(\vecfieldb f))-(\vecfieldb(\vecfielda f))$
for all real valued smooth functions $f$ defined on $\mana$. A \emph{vector
field $\vecfieldcurvea$ along a curve $\curve:\interval\rightarrow\mana$}
is a differentiable mapping that associates to every $t\in\interval$
a tangent vector $\vecfieldcurvea(t)\in T_{\curve(t)}\mana$.
\end{defn}

\begin{defn}[Riemannian manifold]
\label{Definition:Riemannian-metric}A \emph{Riemannian metric} $\inprod{}{}$
or $g$ on a differentiable manifold $\mana$ is a correspondence
which associates to each point $\manpta$ of $\mana$ an inner product
$g_{\manpta}:=\inprod{}{}_{\manpta}$ on a tangent space $T_{\manpta}\mana$,
with $\inprod{}{}_{\manpta}$ varying differentially in the following
sense: if $\charta:U\subset\realset^{n}\rightarrow\mana$ is a system
of coordinates (or chart) around $\manpta$, with $\charta(\parpta_{1},\parpta_{2},...,\parpta_{n})=\manpta\in\charta(U)$
and $\partial/\partial\parpta_{i}(\manpta)=d\charta_{\manpta}(0,...,0,1,0,...0),$
then 
\[
g_{i,j}\left(\parpta_{1},\parpta_{2},...,\parpta_{n}\right)=\inprod{\frac{\partial}{\partial\parpta_{i}}(\manpta)}{\frac{\partial}{\partial\parpta_{j}}(\manpta)}_{\manpta}
\]
is a differentiable function on $U$ \cite{DoCarmo1992}.\emph{ }We
delete the index $\manpta$ in the functions $g_{\manpta}$ and $\inprod{}{}_{\manpta}$
whenever there is no possibility of confusion.

The pair $(\mana,g)$ is called a \emph{Riemannian manifold }\cite{Absil2008}.
For simplicity, we can also denote the Riemannian manifold $(\mana,g)$
by the set $\mana$.
\end{defn}

\begin{defn}[Riemannian gradient \cite{DoCarmo1992}]
 Let $\mana$ be a Riemannian manifold. Given a smooth function $f:\mana\rightarrow\mathbb{R}$,
the \emph{Riemannian gradient} of $f$ at $\man x$, denoted by $\rgrad f(\man x)$
is defined as the unique element of $T_{\man x}\mana$ that satisfies
\[
\left\langle \rgrad f(\man x),v\right\rangle _{\man x}=df_{\man x}(v),\ \forall v\in T_{\man x}\mana.
\]
\end{defn}
\begin{defn}[Critical point \cite{PT:06}]
 Let $\mana$ and $\manb$ be smooth manifolds. If $f:\mana\rightarrow\manb$
is a smooth map, then a point $\man x\in\mana$ is a \emph{critical
point} of $f$ if $df_{\man x}:T_{\man x}\mana\rightarrow T_{\man{f(x)}}\manb$
is not surjective. In the particular case that $\manb=\mathbb{R}$,
then the critical points of $f$ are exactly the points $\man x$
which $df_{\man x}=0$. Moreover, if $\mana$ is a Riemannian manifold,
the critical points are the points $\man x\in\mana$ such that $\rgrad f(\man x)=0$. 
\end{defn}
\begin{defn}[Affine connection \cite{DoCarmo1992}]
\label{thm:Affine-connection}An \emph{affine connection $\nabla$
}on a differentiable manifold $\mana$ is a mapping $\nabla:\vfieldset{\mana}\times\vfieldset{\mana}\rightarrow\vfieldset{\mana}$
which is denoted by $(\vecfielda,\vecfieldb)\mapsto\nabla_{\vecfielda}\vecfieldb$
and which satisfies the following properties, for $\vecfielda$, $\vecfieldb$,
$\vecfieldc$ $\in$ $\vfieldset{\mana}$ and $f$, $g$ $\in$ $\diffset{\mana}$:
\begin{enumerate}
\item $\nabla_{f\vecfielda+g\vecfieldb}\vecfieldc=f\nabla_{\vecfielda}\vecfieldc+g\nabla_{\vecfielda}\vecfieldc$,
\item $\nabla_{\vecfielda}(\vecfieldb+\vecfieldc)=\nabla_{\vecfielda}\vecfieldb+\nabla_{\vecfielda}\vecfieldc$,
\item $\nabla_{\vecfielda}(f\vecfieldb)=f\nabla_{\vecfielda}\vecfieldb+(\vecfielda f)\vecfieldb$.
\end{enumerate}
If \emph{$\nabla$} satisfies the following additional properties:
\begin{enumerate}
\item $\vecfielda\left\langle \vecfieldb,\vecfieldc\right\rangle =\left\langle \nabla_{\vecfielda}\vecfieldb,\vecfieldc\right\rangle +\left\langle \vecfieldb,\nabla_{\vecfielda}\vecfieldc\right\rangle $,
for all $\vecfielda$,$\vecfieldb$, $\vecfieldc\in\vfieldset{\mana}$,
\item $\nabla_{\vecfielda}\vecfieldb-\nabla_{\vecfieldb}\vecfielda=[\vecfielda,\vecfieldb]$,
for all $\vecfielda$,$\vecfieldb\in\vfieldset{\mana}$,
\end{enumerate}
then \emph{$\nabla$} is known as the \emph{Riemannian connection}
of $\mana$. The Levi-Cevita theorem \cite{DoCarmo1992} says that
any Riemannian manifold has a Riemannian connection and it is unique. 
\end{defn}

\begin{thm}[Covariant derivative \cite{DoCarmo1992}]
\label{thm:covariant-derivative}\label{thm:covariant-derivative-1}Let
$\mana$ be a differentiable manifold with an affine connection $\nabla$.
There exists a unique correspondence which associates to a vector
field $\vecfieldcurvea$ along the differentiable curve $\curve:\interval\rightarrow\mana$
another vector field $D\vecfieldcurvea/dt$ along $\curve$, called
the covariant derivative of $\mana$ along $\curve$, such that:
\end{thm}
\begin{enumerate}
\item $\frac{D}{dt}\left(\vecfieldcurvea+\vecfieldcurveb\right)=\frac{D\vecfieldcurvea}{dt}+\frac{D\vecfieldcurveb}{dt}$;
\item $\frac{D}{dt}\left(\vecfieldcurvea V\right)=\frac{df}{dt}\vecfieldcurvea+f\frac{D\vecfieldcurvea}{dt}$,
where $f$ is a differentiable function on $\interval$;
\item if $\vecfieldcurvea$ is induced by a vector field $\vecfieldc\in\vfieldset{\mana}$,
i.e., $\vecfieldcurvea(t)=\vecfieldc(\curve(t))$, then $D\vecfieldcurvea/dt=\nabla_{\curve'(t)}\vecfieldc$.
\end{enumerate}

\begin{defn}[Parallel Transport \cite{DoCarmo1992}]
\label{def:parallel-transport}Let $\mana$ be a differentiable manifold
with an affine connection $\nabla.$ A vector field $\vecfieldcurvea$
along a curve $\curve:\interval\rightarrow\mana$ is called \emph{parallel}
when 
\[
\frac{D\vecfieldcurvea}{dt}(t)=0,\text{ for all }t\in\interval.
\]
Moreover, let $\curve$ be differentiable and $\tanpta_{0}$ a vector
tangent to $\mana$ at $\curve(t_{0})$, $t_{0}\in\interval$. Then
there exists a unique parallel vector field $\vecfieldcurvea$ along
$\curve$, such that $\vecfieldcurvea(t_{0})=\tanpta_{0}$; $V(t)$
is called the \emph{parallel transport} of $\vecfieldcurvea(t_{0})$
along $\curve$.
\end{defn}

\begin{defn}[Geodesic \cite{DoCarmo1992}]
\label{def:geodesics}A parameterized curve $\geod:\interval\rightarrow\mana$
is a \emph{geodesic at} $t_{0}\in\interval$ if 
\[
\frac{D}{dt}(\curve'(t))=0
\]
 at the point $t_{0}$; if $\geod$ is a geodesic at $t$, for all
$t\in\interval$, we say that $\geod$ is a \emph{geodesic} \cite{DoCarmo1992}.
If the definition domain of all geodesics of $\mana$ can be extended
to $\realset$, then $\mana$ is said to be \emph{geodesically-complete}.
\end{defn}

\begin{defn}[Exponential and logarithm mappings \cite{Pennec2006}]
\label{def:exponential-map}Consider a point $\manpta\in\mana$ and
let $\tansubseta\subset T_{\manpta}\mana$ be an open set of $T_{\manpta}\mana$.
For a given vector $v\in V$ and $1\in\interval$, consider the \emph{geodesic}
$\geod:\interval\rightarrow\mana$ passing through $\manpta$ with
initial velocity $\alpha'(0)=v$. Then the mapping $\rexp_{\manpta}:V\rightarrow\mana$
defined by $v\mapsto\alpha(1)$ is well-defined \cite{DoCarmo1992}
and is called the (Riemannian) \emph{exponential mapping} on $\tansubseta$.

The mapping $\rexp_{\manpta}$ is differentiable, and there is a neighborhood
$\submana$ of $\manpta$ such that the exponential map at $\manpta$
is a diffeomorphism from the tangent space to the manifold. For $\submana$
being this neighborhood and $\manpta,\manptb\in\submana$, $\manptb=\rexp_{\manpta}(v)$,
then the inverse mapping $\rlog_{\manpta}:\submana\rightarrow T_{\manpta}\mana$
defined by $\manptb\mapsto\tanpta$ is called the (Riemannian) \emph{logarithm
mapping}. For brevity, we can also write $\rlogb{\manpta\manptb}$
in the place of $\rlog_{\manpta}(\manptb)$.
\end{defn}
\begin{defn}[Riemannian curvature tensor and sectional curvatures \cite{DoCarmo1992}]
 \label{def:Riemannian_curvature} Let $\mathbb{X}(\mana)$ be the
set of mappings from $\vfieldset{\mana}$ to $\vfieldset{\mana}$.
The \emph{Riemannian curvature tensor} $R$ of a differentiable manifold
$\mana$ is the correspondence $R:\vfieldset{\mana}\times\vfieldset{\mana}\rightarrow\mathbb{X}(\mana)$
that associates to each pair of vector fields $\vecfielda,\vecfieldb\in\vfieldset{\mana}$
the application $R(\vecfielda,\vecfieldb):\vfieldset{\mana}\rightarrow\vfieldset{\mana}$
given by 
\[
R(\vecfielda,\vecfieldb)\vecfieldc\coloneqq\nabla_{\vecfieldb}\nabla_{\vecfielda}\vecfieldc-\nabla_{\vecfielda}\nabla_{\vecfieldb}\vecfieldc+\nabla_{[\vecfielda,\vecfieldb]}\vecfieldc,
\]
where $\nabla$ is the Riemannian connection of $\mana$. A notion
closely related to the Riemannian curvature tensor is the sectional
curvatures of $\mana$. Given two linearly independent tangent vectors
$u$ and $v$ at the same point, the expression 
\[
K(u,v)\coloneqq\frac{\left\langle R(u,v)u,v\right\rangle }{\left\langle u,u\right\rangle \left\langle v,v\right\rangle -\left\langle u,v\right\rangle ^{2}}
\]
 does not depend on the choice of $u,v$, but only on the the subspace
$\sigma$ spanned by them \cite{DoCarmo1992}. Given a point $\man p\in\mana$
and a bidimensional subspace $\sigma$ of $T_{\man p}\mana$, the
real number $K(u,v)=K(\sigma)$ where $\{u,v\}$ is any basis of $\sigma$,
is the \emph{sectional curvature} of $\sigma$ in $\man p$.
\end{defn}

\subsection{Proof of Theorem \ref{thm:Euclidean-to-Riemannian-sigma-rep}}

\label{proof:theoremRiSRs}

Suppose $\man{\spa}$ is a Ri$l$th$\nsp\sigma$R of $\man{\rva}$.
Then, from\emph{ }(\ref{eq:Riemannian-sigma-rep-definition-condition}),
(\ref{eq:sigma-rep-def-weights-condition}) is satisfied. Because
$\man{\spa}$ is a Ri$l$th$\nsp\sigma$R of $\man{\rva}$, from (\ref{eq:Riemannian-sigma-rep-definition-mean-condition}),
$\man{\mean{\rva}}$ is a Riemannian sample mean of $\man{\spa}$
and, therefore, from (\ref{eq:Riemannian-sample-mean-definition}),
$\man{\mean{\rva}}$ minimizes the function
\[
g(\man{\state}):=\sum_{i=1}^{\nsp}\weight_{i}^{m}\dist^{2}\left(\man{\state},\rexp_{\mean{\man{\rva}}}\spa_{i}\right).
\]
The function $g\circ\rexp_{\mean{\man{\rva}}}:\maxdefdomain(\mean{\man{\rva}})\subset T_{\mean{\man{\rva}}}\mana\rightarrow[0,\infty)$
is a real valued function defined in a subset of the vector space
$T_{\mean{\man{\rva}}}\mana$. Since $\maxdefdomain(\mean{\man{\rva}})$
is convex by hypothesis and its second derivative is positive, then
$g\circ\rexp_{\mean{\man{\rva}}}$ is a strictly convex function.
Because it is also a differentiable function, $g\circ\rexp_{\manpta}$
has an unique minimum $\state^{*}\in\maxdefdomain(\manpta)$ and it
is a critical point of $g\circ\rexp_{\manpta}$. Thus
\begin{equation}
[0]_{n\times1}=\left.\frac{d\big(g\circ\rexp_{\mean{\man{\rva}}}\big)(x)}{dx}\right|_{x=\state^{*}}\Leftrightarrow\state^{*}=\sum_{i=1}^{\nsp}\weight_{i}^{m}\spa_{i}.\label{eq:sample-mean-in-tangent-space}
\end{equation}
By Theorem~7.9 of \cite{NB:13} , $\man{\mean{\rva}}$ is the unique
minimum and critical point of $g$. Thus $\rlog_{\manpta}\man{\mean{\rva}}$
is a critical point of $g\circ\rexp_{\manpta}$. and
\begin{equation}
[0]_{n\times1}=\rlogb{\mean{\man{\rva}}\man{\mean{\rva}}}=\state^{*}=\sum_{i=1}^{\nsp}\weight_{i}^{m}\spa_{i}=:\smean_{\spa}.\label{eq:sample-mean-in-tangent-space-1}
\end{equation}
Hence, (\ref{eq:sigma-rep-def-mean-condition}) is satisfied.

Now let us prove the converse for the mean. Suppose all points $\spa_{i}$
belong to the domain of $\rexp_{\mean{\man{\rva}}}$, and that $\spa$
is an $l$th$\nsp\sigma$R of $\rva:=\rlogb{\mean{\man{\rva}}\man{\rva}}$.
Define the set
\begin{multline}
\man{\spa}:=\{\rexp_{\mean{\man{\rva}}}\spa_{i},\weight_{i}^{m},\weight_{i}^{c,j},\weight_{i}^{cc,j}|\man{\spa}_{i}\in\mana;\\
\,\weight_{i}^{cc,j},\weight_{i}^{cc,j},\weight_{i}^{cc,j}>0\}{}_{i=1}^{\nsp}\label{eq:sigma-rep-proof}
\end{multline}
Then, from (\ref{eq:sigma-rep-def-weights-condition}) and (\ref{eq:sigma-rep-proof})\emph{,
}(\ref{eq:Riemannian-sigma-rep-definition-condition}) is satisfied.
From (\ref{eq:sample-mean-in-tangent-space}) and (\ref{eq:sample-mean-in-tangent-space-1}),
we have that $\rexp_{\mean{\man{\rva}}}(\smean_{\spa})=\rexp_{\mean{\man{\rva}}}(\rlogb{\mean{\man{\rva}}\man{\mean{\rva}}})=\mean{\man{\rva}}$
minimizes $g$ and (\ref{eq:Riemannian-sigma-rep-definition-mean-condition})
is satisfied.

For even $j$, we have, from (\ref{eq:sigma-rep-def-mean-condition})
and (\ref{eq:Riemannian-sample-moment-definition}),
\[
\man{\smoment}_{\man{\spa}}^{j}=\sum_{i=1}^{\nsp}\weight_{i}^{cc,j}\big[(\spa_{i}-\smean_{\spa})\left(\diamond\right)^{T}\big]^{\otimes\frac{j}{2}}=:\mathcal{\smoment}_{\spa}^{j};
\]
and from (\ref{eq:sigma-rep-def-moments-condition}), it follows $\man{\smoment}_{\man{\spa}}^{j}=\mathcal{\smoment}_{\spa}^{j}=\moment_{\rva}^{j};$
for odd $j$, the reasoning is similar. The remaining is straightforward.

\subsection{Proof of Corollary \ref{cor:Riemannian-sr-minimum-numbers}}

\label{proof:Corollary_particularRiSRs}

From Theorem \ref{thm:Euclidean-to-Riemannian-sigma-rep}, $\spa$
is a normalized\emph{ }$l$th$\nsp\sigma$R of\emph{ }$\rva\sim(\rlog_{\mean{\man{\rva}}}(\man{\rva}),\man{\cov}_{\man{\rva}\man{\rva}})^{n}$.
From Corollary 1 of \cite{Menegaz2015}, it follows that i) $\nsp\geq r+1$;
and ii), if $\man{\spa}$ is symmetric, then $\nsp=2r$. The remaining
of the proof is a direct consequence of Theorem \ref{thm:Euclidean-to-Riemannian-sigma-rep}.

\subsection{Proof of Proposition \ref{prop:Addition of Riemannian random points}}

\label{proof:proposition1}

A Riemannian mean $\mean{\man{\rva}}$ of $\man{\rva}$ is such that
it solves (\ref{eq:Riemannian-mean-definition}). Consider the following
optimization problem
\begin{align}
\mbox{minimize } & \tilde{g}(\tilde{c}):=g\circ\rexp_{\bar{\manpta}}\left(\rlogb{\mean{\manpta}\manptcentral}\right)=\sd_{\rlogb{\mean{\manpta}\manpta}+p}^{2}(\tilde{c})\nonumber \\
\mbox{subject to } & \manptcentral\in\mana;\label{eq:Riemannian-addition-optimization-problem2}
\end{align}
From a reasoning similar to the proof of Theorem \ref{thm:Euclidean-to-Riemannian-sigma-rep}
(Appendix \ref{proof:theoremRiSRs}), if $\tilde{c}$ solves (\ref{eq:Riemannian-addition-optimization-problem2}),
then $\rlog_{\mean{\manpta}}^{-1}\tilde{c}=\rexp_{\bar{\manpta}}\tilde{c}$
solves (\ref{eq:Riemannian-mean-definition}), and $\mean{\man{\rva}}=\rexp_{\bar{\manpta}}\tilde{c}$.
Since $\sd_{\rlogb{\mean{\manpta}\manpta}+p}^{2}(\tilde{c})$ is the
variance of $\rlogb{\mean{\manpta}\manpta}+p$ it follows that $\ev_{\rlogb{\mean{\manpta}\manpta}+p}\{\rlogb{\mean{\manpta}\manpta}+p\}=\mean p$
minimizes $\tilde{g}\left(\manptcentral\right)$; thus $\mean{\man{\rva}}:=\rexp_{\bar{\manpta}}\mean p$.
For the covariance part, we have
\[
\man{\cov}_{\man{\rva}\man{\rva}}:=\int_{\mana-\cutlocus(\mean{\man{\rva}})}\rlogb{\mean{\man{\rva}}\man x}\big(\rlogb{\mean{\man{\rva}}\man x}\big)^{T}\man{\pdf}_{\man{\rva}}(\man x)d\mana\left(\man x\right)=\man{\cov}_{\manpta\manpta}+\cov_{pp}.
\]

\subsection{Proof of Theorem \ref{thm:General-to-Hauberg-UKF}}

\label{proof:TheoremRiUKFcorrection}

First, by considering $\man c_{\man{\state}}=\man b_{\man{\state}}=\mean{\man{\state}}_{k|k-1}$
and $\man c_{\man{\meas}}=\man b_{\man{\meas}}=\mean{\man{\meas}}_{k|k-1}$
in the definitions of $\man{\cov}_{\man{\state}\man{\state},**}^{k|k-1}$,
$\man{\cov}_{\man{\meas}\man{\meas},**}^{k|k-1}$, and $\man{\cov}_{\man{\state}\man{\meas},**}^{k|k-1}$,
(\ref{eq:UKF-Riemannian-CartesianSolution-3}) yields
\begin{equation}
\man{\kgain}_{k,**}=\left[\begin{array}{cc}
[0]_{\dimstate\times\dimstate} & \man{\kgain}_{k}\\{}
[0]_{\dimmeas\times\dimstate} & [0]_{\dimmeas\times\dimmeas}
\end{array}\right],\label{eq:proof-general-to-Hauberg-UKF-3}
\end{equation}
and substituting $\man c_{\man{\state}}=\man b_{\man{\state}}=\est{\man{\state}}_{k|k-1}$,
$\man c_{\man{\meas}}=\man b_{\man{\meas}}=\est{\man{\meas}}_{k|k-1}$,
and (\ref{eq:proof-general-to-Hauberg-UKF-3}) into (\ref{eq:UKF-Riemannian-CartesianSolution-4})
gives $\mean{\state}_{k|k,**}^{T_{\man c}\mana_{\man{\state},\man{\meas}}}=\big[\begin{array}{c}
\man{\kgain}_{k}\rlog_{\mean{\man{\meas}}_{k|k-1}}(\outcome{\man{\meas}}_{k}),[0]_{\dimmeas\times1}\end{array}\big]^{T};$ consequently, from (\ref{eq:augmented-corrected-state-to-non-augmented}),
$x_{k|k}^{TM}:=[x_{k|k,**}^{T_{\man c}\mana_{\man{\state},\man{\meas}}}]_{1:\dimstate,1}=\man{\kgain}_{k}\rlog_{\man{\meas}_{k|k-1}}\left(\man{\meas}_{k}\right)$.
Second, considering (\ref{eq:proof-general-to-Hauberg-UKF-3}) into
(\ref{eq:UKF-Riemannian-CartesianSolution-6}) yields
\[
\man{\cov}_{\man{\state}\man{\state},**}^{k|k,T_{\man c}M}=\diag\Big(\man{\cov}_{\man{\state}\man{\state}}^{k|k-1}-\man{\kgain}_{k}\big(\man{\cov}_{\man{\meas}\man{\meas}}^{k|k-1}\big)^{-1}\man{\kgain}_{k}^{T},[0]_{\dimmeas\times\dimmeas}\Big);
\]
and, from (\ref{eq:augmented-corrected-state-to-non-augmented}),
it follows that $\man{\cov}_{\man{\state}\man{\state},**}^{k|k,\est{\man{\state}}_{k|k-1}}=\man{\cov}_{\man{\state}\man{\state}}^{k|k-1}-\man{\kgain}_{k}(\man{\cov}_{\man{\meas}\man{\meas}}^{k|k-1})^{-1}\man{\kgain}_{k}^{T}$.

\subsection{Notation and Acronyms}

\label{subsec:Notation-and-Acronyms}

Throughout this paper, we use the following notations:
\begin{itemize}
\item for a matrix $A$, $(A)\left(\diamond\right)^{T}$ stands for $(A)\left(A\right)^{T}$,
and $\sqrt{A}$ for a square-root matrix of $A$ such that $A=\sqrt{A}\sqrt{A}^{T}$.
\item $\otimes$ stands for the Kronecker product operator, and $A^{\otimes n}:=A\otimes\cdots\otimes A$. 
\item {\small{}$[A]_{p\times q}$} stands for a block matrix consisting
of the matrix $A$ being repeated $p$ times in the rows and $q$
times in the columns.
\item {\small{}$[A]_{i_{1}:i_{2},j_{1}:j_{2}}$} stands for a sub-matrix
of the matrix $A$ formed by the rows $i_{1}$ to $i_{2}$ and the
columns $j_{1}$ to $j_{2}$ of $A$.
\item $\interval$ stands for an open interval in $\realset$.
\end{itemize}

Below, we provide a list of acronyms and parts of acronyms—these parts
end with an '-' and are followed by examples—along with their meaning.
There are other acronyms in the text that can be composed by i) concatenating
some items below (e.g., Mi$\sigma$R {[}Mi- with $\sigma$R{]} standing
for Minimum $\sigma$Representation) , and ii) adding Ri- (standing
for Riemannian; e.g., RiMi$\sigma$R\ {[}Ri- with Mi$\sigma$R{]}
standing for RiMi$\sigma$R): 

\begin{itemize}
\item \textbf{AdUKF}: Additive Unscented Kalman Filter.
\item \textbf{AuUKF}: Augmented Unscented Kalman Filter.
\item \textbf{EKF}: Extended Kalman Filter.
\item \textbf{HoMiSy-}: Homogeneous Minimum Symmetric- (e.g., HoMiSy$\sigma$R,
RiHoMiSy$\sigma$R, RiHoMiSyAdUKF, RiHoMiSyAuUKF).
\item \textbf{KF}: Kalman Filter.
\item \textbf{$l$th$N\sigma$R}: $l$th order $N$ points $\sigma$-representation.
\item \textbf{$l$UT}: $l$th order UT.
\item \textbf{Mi-}: Minimum- (e.g., RiMi$\sigma$R, RiMiAdUKF, RiMiAuUKF).
\item \textbf{MiSy-}: Minimum Symmetric- (e.g., RiMiSy$\sigma$R, RiMiSyAdUKF,
RiMiSyAuUKF).
\item \textbf{RhoMi-}: Rho Minimum- (e.g., RiRhoMi$\sigma$R, RiRhoMiAdUKF,
RiRhoMiAuUKF).
\item \textbf{$\sigma$R}: $\sigma$-Representation.
\item \textbf{UKF}: Unscented Kalman Filter.
\item \textbf{UT}: Unscented Transformation.
\end{itemize}
%
\begin{IEEEbiography}[{\includegraphics[width=1in,height=1.25in,clip,keepaspectratio]{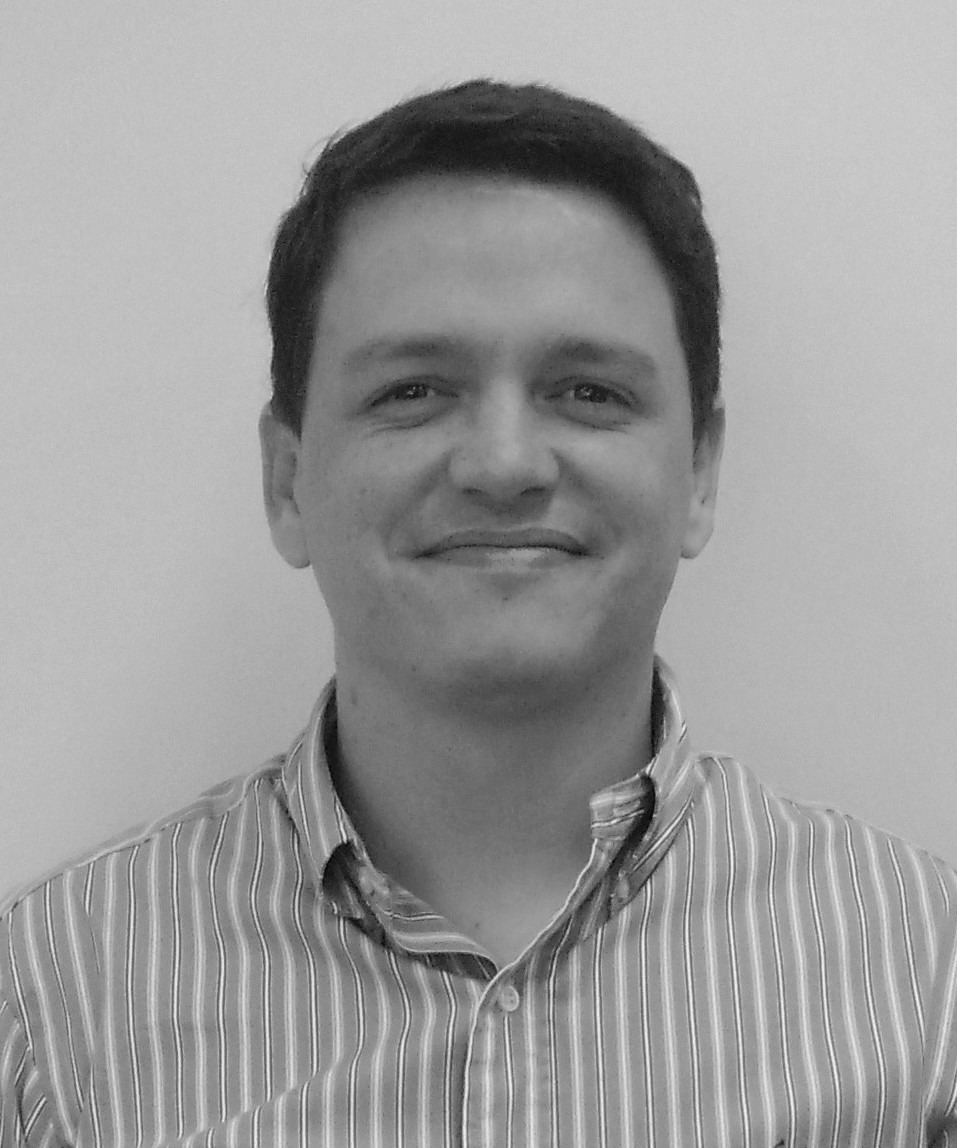}}]{Henrique M. T. Menegaz} received the B.S. degree in Electrical Engineering from the Universidade de Brasília (UnB), Brazil, in 2007. He received the M.S. and Ph.D. degrees in Engineering of Electronic Systems and Automation from the UnB in 2011 and 2016 respectively. He is currently an Assistant Professor with Faculdade Gama, UnB. His major field of study is filtering of nonlinear dynamic systems and their applications. 
\end{IEEEbiography}
\begin{IEEEbiography}[{\includegraphics[width=1in,height=1.25in,clip,keepaspectratio]{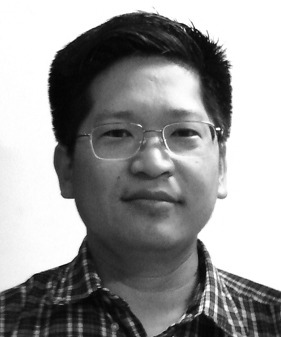}}]{João Y. Ishihara} received the Ph.D. degree in Electrical Engineering from the University of São Paulo, Brazil, in 1998.
He is currently an Associate Professor at the University of Brasília, Brazil. His research interests include robust filtering and control theory, singular systems, and robotics.
\end{IEEEbiography}
\begin{IEEEbiography}[{\includegraphics[width=1in,height=1.25in,clip,keepaspectratio]{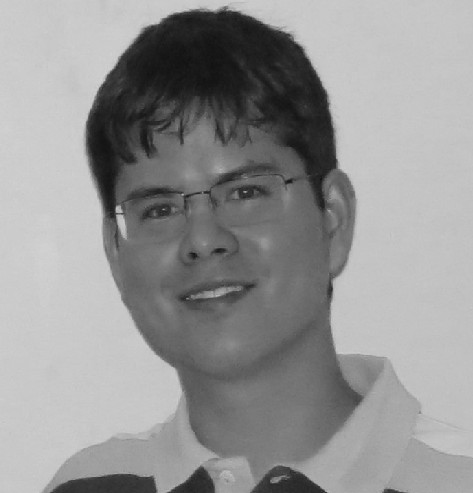}}]{Hugo T. M. Kussaba} received the B.S. degree in control engineering and the M.S. degree in engineering of electronic systems and automation from the University of Brasília (UnB), Brazil, in 2012 and 2014 respectively. Currently, he is a  Ph.D. student at the same university. His research interests include control and estimation on Lie groups, robust control and linear matrix inequalities, and hybrid dynamical systems.
\end{IEEEbiography}

\end{document}